\documentclass[sn-mathphys-num]{sn-jnl}% Math and Physical Sciences Numbered Reference Style 
\usepackage{graphicx}%
\usepackage{multirow}%
\usepackage{amsmath,amssymb,amsfonts}%
\usepackage{amsthm}%
\usepackage{mathrsfs}%
\usepackage[title]{appendix}%
\usepackage[dvipsnames]{xcolor}%
\usepackage{textcomp}%
\usepackage{manyfoot}%
\usepackage{booktabs}%
\usepackage{algorithm}%
\usepackage{algorithmicx}%
\usepackage{algpseudocode}%
\usepackage{listings}%
\usepackage{subcaption}
\usepackage{stmaryrd}
\usepackage{bm}
\usepackage{empheq}
\newcommand{\R}{\mathbb{R}}
\newcommand{\M}{\mathcal{M}}
\newcommand{\V}{\mathcal{V}}
\newcommand{\I}{\mathcal{I}}
\newcommand{\Po}{\bm{\mathcal P}}
\newcommand{\Ro}{\bm{\mathcal R}}

\newcommand*\diff{\mathop{}\!\mathrm{d}}
\newcommand{\dott}[2]{\left\langle#1,#2\right\rangle}
\newcommand{\norm}[1]{\left\lVert #1\right\rVert}
\newcommand{\intervalle}[4]{\mathopen{#1}#2\mathclose{},#3\mathclose{#4}}
\newcommand{\ff}[2]{\intervalle{[}{#1}{#2}{]}}
\newcommand{\oo}[2]{\intervalle{(}{#1}{#2}{)}}

\floatname{algorithm}{Algorithm}

\algnewcommand{\algorithmicendif}{\textbf{end}}
\algblockdefx[IF]{If}{EndIf}[1]{\algorithmicif\ #1\ \algorithmicthen}{\algorithmicendif}
\algblockdefx[For]{for}{EndFor}[1]{\algorithmicif\ #1\ \algorithmicthen}{\algorithmicendif}

\newcommand{\review}[1]{{\textcolor{black}{#1}}}
\newcommand{\modif}[1]{{\textcolor{black}{#1}}}
\newcommand{\typo}[1]{{\textcolor{black}{#1}}}

\theoremstyle{thmstyleone}%
\newtheorem{theorem}{Theorem}%  meant for continuous numbers
\newtheorem{proposition}{Proposition}
\newtheorem{assumption}{Assumption}
\newtheorem{lemma}{Lemma}

\theoremstyle{thmstyletwo}%
\newtheorem{remark}{Remark}%

\theoremstyle{thmstylethree}%

\begin{document}

\title[Article Title]{\modif{High-order BUG dynamical low-rank integrators based on explicit Runge--Kutta methods}}


\author[1]{\fnm{Fabio} \sur{Nobile}}
\author*[1]{\fnm{S\'ebastien} \sur{Riffaud}}\email{sebastien.riffaud@epfl.ch}

\affil[1]{\orgdiv{CSQI Chair}, \orgname{\'Ecole Polytechnique F\'ed\'erale de Lausanne}, \city{1015 Lausanne}, \country{Switzerland}}


\abstract{In this work, we \typo{introduce} high-order Basis-Update \& Galerkin (BUG) integrators based on explicit Runge--Kutta methods for large-scale matrix differential equations. These dynamical low-rank integrators \typo{extend} the BUG integrator~\cite{ceruti2022rank} \typo{to arbitrary explicit Runge--Kutta schemes} by performing a BUG step at each stage of the method. The resulting Runge--Kutta BUG \typo{(RK--BUG)} integrators are robust \typo{with respect to} small singular values, \typo{fully forward in time}, \modif{and high-order accurate, while enabling conservation and rank adaptivity}. We \typo{prove} that \typo{RK--BUG} integrators retain the order of convergence of the \typo{underlying} Runge--Kutta method until the error reaches a plateau corresponding to the low-rank truncation error, which vanishes as the rank becomes full. \typo{This theoretical analysis is supported by several numerical experiments}. The results demonstrate the high-order convergence of the \typo{RK--BUG} integrator and its superior accuracy compared to other \typo{existing} dynamical low-rank integrators.}




\maketitle

\section{Introduction}
\label{sec:1}

Dynamical low-rank approximations (DLRAs) enable a significant reduction in the computational cost associated with \typo{the} numerical integration of large-scale matrix differential equations:
\begin{equation}
\label{eq:ode}
\typo{\dot{\mathbf X}} = {\mathbf F}(t,\typo{{\mathbf X}}), \qquad \typo{{\mathbf X}}(0) = \typo{{\mathbf X}}_0 \in \R^{n \times m},
\end{equation}
which \typo{arise} in many applications, such as kinetic equations~\cite{bernard2018reduced,einkemmer2018low,einkemmer2020low,einkemmer2021mass,coughlin2022efficient,einkemmer2024accelerating,einkemmer2024review} due to the large phase space, stochastic simulations~\cite{sapsis2009dynamically,babaee2017robust,musharbash2018dual,feppon2018dynamically,musharbash2020symplectic,patil2020real,kazashi2021existence,kazashi2021stability,kazashi2025dynamical} due to the repeated evaluation of the solution for different realizations of the random terms, or sequential parameter estimation~\cite{kressner2011low,weinhandl2018linear,benner2024low,riffaud2024low} when tracking solutions for several parameter values. The main idea of DLRAs \typo{is to approximate} the solution \typo{${\mathbf X}$} by a low-rank matrix ${\mathbf Y}$ in an SVD-like \typo{factorization}:
\begin{equation}
\label{eq:sol}
{\mathbf Y}(t) = {\mathbf U}(t)\,{\mathbf S}(t)\,{\mathbf V}(t)^\top \in \R^{n \times m},
\end{equation}
where ${\mathbf U} \in \R^{n \times r}$ and ${\mathbf V} \in \R^{m \times r}$ are orthonormal matrices, ${\mathbf S} \in \R^{r \times r}$ is an \typo{invertible square} matrix (not necessarily diagonal), and $r \leq \min\{n,m\}$ \typo{denotes} the rank of ${\mathbf Y}$.

A challenging issue concerns the time-integration of the \typo{low-rank factorization~\eqref{eq:sol}}. Perhaps the most intuitive approach is to consider the \typo{system of} differential equations, derived in~\cite{koch2007dynamical}, that describes the evolution of the low-rank \typo{factors $({\mathbf U}, {\mathbf S}, {\mathbf V})$} over time:
\begin{equation}
\label{eq:low-rank_ode}
\left\{
\begin{aligned}
\dot{\mathbf U} &= \left({\mathbf I} - {\mathbf U}{\mathbf U}^\top \right) {\mathbf F}(t,{\mathbf Y}){\mathbf V}{\mathbf S}^{-1},\\
\dot{\mathbf S} &= {\mathbf U}^\top {\mathbf F}(t,{\mathbf Y}){\mathbf V},\\
\dot{\mathbf V} &= \left({\mathbf I} - {\mathbf V}{\mathbf V}^\top \right){\mathbf F}(t,{\mathbf Y})^\top{\mathbf U}{\mathbf S}^{-T},
\end{aligned}
\right.
\end{equation}
where \typo{the orthogonality constraints ${\mathbf U}^\top\dot{\mathbf U} = {\mathbf 0}$ and ${\mathbf V}^\top\dot{\mathbf V} = {\mathbf 0}$ have been imposed for uniqueness}. Unfortunately, this system involves the inverse of ${\mathbf S}$, which can \typo{lead to} severe restrictions \typo{on the time-step when} the singular values of ${\mathbf S}$ are small, since the step size of standard time-integration schemes must be proportional to the smallest nonzero singular value. An equivalent formulation \typo{of}~\eqref{eq:low-rank_ode} is obtained by projecting the matrix differential equation~\eqref{eq:ode} onto the tangent space of the manifold $\M_r$ of rank-$r$ matrices:
\begin{equation}
\label{eq:projected_ode}
\dot{\mathbf Y} = \Po_{\mathbf Y}\big({\mathbf F}(t,{\mathbf Y})\big), \qquad {\mathbf Y}(0) = {\mathbf Y}_0 \in \M_r,
\end{equation}
where 
\begin{equation}
\label{eq:tangent-space_projection}
\Po_{\mathbf Y}(\typo{{\mathbf Z}}) = {\mathbf U}{\mathbf U}^\top\typo{{\mathbf Z}} - {\mathbf U}{\mathbf U}^\top\typo{{\mathbf Z}}{\mathbf V}{\mathbf V}^\top + \typo{{\mathbf Z}}{\mathbf V}{\mathbf V}^\top
\end{equation}
stands for the orthogonal projection of $\typo{{\mathbf Z}} \in \R^{n \times m}$ onto the tangent space at ${\mathbf Y}={\mathbf U}{\mathbf S}{\mathbf V}^\top \in \M_r$.

In the \textit{\typo{projector}-splitting integrator}~\cite{lubich2014projector,kieri2016discretized}, \typo{the} tangent-space projection~\eqref{eq:tangent-space_projection} is \typo{decomposed} into three \typo{substeps (associated with the updates of ${\mathbf U}$, ${\mathbf S}$, and ${\mathbf V}$) which are integrated sequentially} using the Lie--Trotter or Strang splitting. The resulting \typo{dynamical} low-rank integrator is robust \typo{with respect to} small singular values, but the \typo{S-step} integrates the solution backward in time, which \typo{may cause} instabilities for parabolic and hyperbolic problems~\cite{kusch2023stability}. \modif{To avoid backward time-integration, the \textit{unconventional integrator}~\cite{ceruti2022unconventional} first updates ${\mathbf U}$ and ${\mathbf V}$ in parallel using the tangent-space projection, and then updates ${\mathbf S}$ through forward time-integration via a Galerkin projection. This approach removes the backward time-integration and remains robust with respect to small singular values, but is limited to a fixed rank~$r$. The \textit{Basis-Update \& Galerkin (BUG) integrator}~\cite{ceruti2022rank,ceruti2024parallel} is a rank-adaptive variant of the unconventional integrator. In this approach, ${\mathbf U}$ and ${\mathbf V}$ are first updated and augmented using the tangent-space projection, then the augmented matrix $\widehat{\mathbf S}$ is updated via a Galerkin projection, and finally the augmented rank is truncated. As a result, the BUG integrator is robust with respect to small singular values, fully forward in time, and enables rank adaptivity, but its convergence order is limited to one due to the first-order Lie--Trotter splitting.}

\modif{In recent years, several high-order dynamical low-rank integrators that are robust, forward-only, and rank-adaptive have been proposed:
\begin{itemize}
\item the \textit{midpoint BUG integrator}~\cite{ceruti2024robust,kusch2024second}, which is a second-order extension of the BUG integrator based on the midpoint rule;
\item the \textit{projected Runge--Kutta (PRK) methods}~\cite{kieri2019projection,carrel2023projected}, which integrate the projected differential equation~\eqref{eq:projected_ode} using high-order explicit Runge--Kutta schemes and include a truncation step that maintains low-rank approximations.
\end{itemize}
Compared to the midpoint BUG integrator, PRK methods achieve higher-order convergence for time-explicit discretizations. However, the Galerkin projection used in the BUG integrator is more accurate than the tangent-space projection employed in PRK methods for approximating ${\mathbf F}$ at the discrete level (see Proposition~\ref{theo:P3})}. \review{Moreover, the BUG integrator can be adapted, with only minor modifications~\cite{einkemmer2023conservation}, to preserve important conservation properties, such as the conservation of mass, momentum, and energy in kinetic equations}.

In this work, we \typo{introduce} high-order BUG integrators \typo{that extend the first-order BUG integrator to arbitrary} explicit Runge--Kutta methods. These \typo{Runge--Kutta BUG (RK--BUG)} integrators are \typo{obtained} by performing a BUG step at each stage of the Runge--Kutta \typo{scheme}. \typo{As a result, RK--BUG integrators} are robust, \typo{forward-only in time}, \modif{and high-order accurate, while enabling rank adaptivity and the preservation of conservation properties}. \typo{In particular}, we prove in Theorem~\ref{theo:T2} that \typo{RK--BUG} integrators retain the order of convergence of the \typo{underlying} Runge--Kutta method until the error reaches a plateau corresponding to the low-rank truncation error, which vanishes as the rank becomes full. \typo{Moreover}, this property holds for any explicit Runge--Kutta \typo{scheme}, allowing in practice the construction of a \typo{wide class} of high-order dynamical low-rank integrators.

The remainder of the paper is organized as follows. Section~\ref{sec:2} \modif{introduces the RK--BUG integrators, starting from the first-order BUG formulation and extending it to high-order explicit Runge--Kutta schemes}. Then, Section~\ref{sec:3} \typo{demonstrates the high-order} convergence of the proposed \typo{RK--BUG} integrators. In Section~\ref{sec:4}, \typo{this theoretical result} is validated \typo{through several numerical experiments}. Finally, Section~\ref{sec:5} draws some conclusions and perspectives.

\section{\modif{RK--BUG integrators}}
\label{sec:2}

\modif{In this section, we introduce the Runge--Kutta Basis-Update \& Galerkin (RK--BUG) integrators, which generalize the first-order BUG integrator to arbitrary explicit Runge--Kutta methods. We begin by recalling its formulation in the time-explicit case and reinterpreting it in a way that naturally extends to higher-order schemes.}

\subsection{\modif{Reinterpretation in the time-explicit case}}

\typo{For the convenience of the reader, we recall here the formulation of the BUG integrator in the time-explicit case, which is based on the forward Euler method.} Let the time be discretized using a fixed step size $h>0$. \typo{A BUG step for integrating} the rank-$r$ solution ${\mathbf Y}_k = {\mathbf U}_k {\mathbf S}_k {\mathbf V}_k^\top$ from time $t_k$ to $t_k + h$ reads:
\begin{enumerate}
\item \textbf{K-step:} Assemble
\begin{equation}
\label{eq:k_step}
{\mathbf K} = \bigl[\,{\mathbf U}_k \;\; {\mathbf F}(t_k,{\mathbf Y}_k){\mathbf V}_k \,\bigr] \in \R^{n \times 2r},
\end{equation}
and compute $\widehat{\mathbf U}_{k+1} \in \R^{n \times \hat r}$, with $\hat r \in \{r,\ldots,\min\{n,m,2r\}\}$,  
as an orthonormal basis of the range of ${\mathbf K}$ (e.g., by QR decomposition),  
in short $\widehat{\mathbf U}_{k+1} = \operatorname{orth}({\mathbf K})$.
\item \textbf{L-step:} Assemble
\begin{equation}
\label{eq:l_step}
{\mathbf L} = \bigl[\,{\mathbf V}_k \;\; {\mathbf F}(t_k,{\mathbf Y}_k)^\top{\mathbf U}_k \,\bigr] \in \R^{m \times 2r},
\end{equation}
and compute the augmented basis $\widehat{\mathbf V}_{k+1} = \operatorname{orth}({\mathbf L}) \in \R^{m \times \hat r}$.
\item \textbf{S-step:} Set
\begin{equation}
\label{eq:s_step}
\widehat{\mathbf S}_{k+1}
= \widehat{\mathbf U}_{k+1}^\top \bigl({\mathbf Y}_k + h\,{\mathbf F}(t_k,{\mathbf Y}_k)\bigr) \widehat{\mathbf V}_{k+1}
\in \R^{\hat r \times \hat r}.
\end{equation}
\item \textbf{Truncation step:}  
\typo{Compute the $r$-truncated singular value decomposition (SVD)
\(
\widehat{\mathbf S}_{k+1} \approx {\boldsymbol\Phi}\,{\boldsymbol\Sigma}\,{\boldsymbol\Psi}^\top,
\)}
where ${\boldsymbol\Phi},{\boldsymbol\Psi} \in \R^{\hat r \times r}$ are orthonormal matrices and ${\boldsymbol\Sigma} \in \R^{r \times r}$ is a diagonal matrix with non-negative \typo{entries} on the diagonal. \typo{Then,} set
\begin{equation}
\label{eq:t_step}
{\mathbf U}_{k+1} = \widehat{\mathbf U}_{k+1} {\boldsymbol\Phi}, \qquad {\mathbf S}_{k+1} = {\boldsymbol\Sigma}, \qquad {\mathbf V}_{k+1} = \widehat{\mathbf V}_{k+1} {\boldsymbol\Psi}.
\end{equation}
\end{enumerate}
\typo{This procedure} is summarized in Algorithm~\ref{alg:BUG}. \typo{For clarity, we denote by $\Ro_{\M_r}({\mathbf Z})$} the projection \typo{(or retraction)} of $\typo{{\mathbf Z}}$ onto \typo{the manifold} $\M_r$ \typo{of rank-$r$ matrices}, which is given by \typo{its} $r$-truncated SVD.

\modif{
A BUG step can be interpreted as follows. The augmented bases $(\widehat{\mathbf U}_{k+1}, \widehat{\mathbf V}_{k+1})$ are first constructed to provide an exact representation of the projected discrete solution of~\eqref{eq:projected_ode}:
\begin{equation*}
\begin{aligned}
\widetilde{\mathbf Z}_{k+1}
  =&\; {\mathbf Y}_k + h\, \Po_{{\mathbf Y}_k}\!\bigl({\mathbf F}(t_k,{\mathbf Y}_k)\bigr) \\
  =&\; {\mathbf U}_k {\mathbf S}_k {\mathbf V}_k^{\top}
  + h\bigl(
     {\mathbf U}_k {\mathbf U}_k^{\top} {\mathbf F}(t_k,{\mathbf Y}_k)
     - {\mathbf U}_k {\mathbf U}_k^{\top} {\mathbf F}(t_k,{\mathbf Y}_k) {\mathbf V}_k {\mathbf V}_k^{\top}
     + {\mathbf F}(t_k,{\mathbf Y}_k) {\mathbf V}_k {\mathbf V}_k^{\top}
     \bigr) \\
    =&\; 
    \begin{bmatrix}
    {\mathbf U}_k & {\mathbf F}(t_k,{\mathbf Y}_k){\mathbf V}_k
    \end{bmatrix}
    \begin{bmatrix}
    {\mathbf S}_k - h\, {\mathbf U}_k^{\top} {\mathbf F}(t_k,{\mathbf Y}_k) {\mathbf V}_k & h\, {\mathbf I} \\
    h\, {\mathbf I} & {\mathbf 0}
    \end{bmatrix}
    \begin{bmatrix}
    {\mathbf V}_k & {\mathbf F}(t_k,{\mathbf Y}_k)^\top {\mathbf U}_k
    \end{bmatrix}^\top.
\end{aligned}
\end{equation*}
Specifically, $\widehat{\mathbf U}_{k+1}$ and $\widehat{\mathbf V}_{k+1}$ are chosen such that they span the column and row spaces of $\widetilde{\mathbf Z}_{k+1}$, that is,
\begin{equation*}
  \widehat{\mathbf U}_{k+1}
  = \operatorname{orth}\bigl([\,{\mathbf U}_k,\;{\mathbf F}(t_k,{\mathbf Y}_k){\mathbf V}_k\,]\bigr), \qquad
  \widehat{\mathbf V}_{k+1}
  = \operatorname{orth}\bigl([\,{\mathbf V}_k,\;{\mathbf F}(t_k,{\mathbf Y}_k)^{\top}{\mathbf U}_k\,]\bigr).
\end{equation*}
Then, ${\widehat{\mathbf S}}_{k+1}$ is computed to best approximate the (unprojected) discrete solution of~\eqref{eq:ode}:
\begin{equation*}
  {\mathbf Z}_{k+1} = {\mathbf Y}_k + h\,{\mathbf F}(t_k,{\mathbf Y}_k), \qquad
  {\widehat{\mathbf S}}_{k+1}
  = \arg\min_{{\mathbf S} \in \R^{\hat r \times \hat r}}\;
    \bigl\|\,\widehat{\mathbf U}_{k+1} \mathbf S \widehat{\mathbf V}_{k+1}^\top
      - {\mathbf Z}_{k+1}\,\bigr\|_F,
\end{equation*}
or, equivalently, by Galerkin projection,
\begin{equation*}
  {\widehat{\mathbf S}}_{k+1}
  = \widehat{\mathbf U}_{k+1}^{\top}
    ({\mathbf Y}_k + h\,{\mathbf F}(t_k,{\mathbf Y}_k))
    \widehat{\mathbf V}_{k+1}.
\end{equation*}
The augmented solution $\widehat{\mathbf Y}_{k+1}:=\widehat{\mathbf U}_{k+1} \widehat{\mathbf S}_{k+1} \widehat{\mathbf V}_{k+1}^\top$ is finally projected onto $\M_r$ via its $r$-truncated SVD, which provides the best rank-$r$ approximation. Compared to the original formulation, this reinterpretation does not rely on any (first-order) splitting of the tangent-space projection, which makes its extension to higher-order explicit Runge--Kutta methods both natural and straightforward.
}

\begin{algorithm}
\caption{BUG integrator based on the forward Euler method (\typo{following}~\cite{ceruti2022rank})}
\label{alg:BUG}
\begin{algorithmic}[1]
\Require {${\mathbf Y}_0$.}
\Ensure {${\mathbf Y}_{1},\ldots,{\mathbf Y}_{N}$.}
\For {$k=0,\ldots,N-1$}
\State ${\mathbf Y}_k := {\mathbf U}_{k} {\mathbf S}_{k} {\mathbf V}_{k}^\top$;
\State ${\mathbf F}_k \gets {\mathbf F}(t_k,{\mathbf Y}_k)$;
\State $\widehat{\mathbf U}_{k+1} \gets \operatorname{orth}\bigl([\,{\mathbf U}_k,\;{\mathbf F}_k{\mathbf V}_k\,]\bigr)$; \Comment{K-step \eqref{eq:k_step}}
\State $\widehat{\mathbf V}_{k+1} \gets \operatorname{orth}\bigl([\,{\mathbf V}_k,\;{\mathbf F}_k^{\top}{\mathbf U}_k\,]\bigr)$; \Comment{L-step \eqref{eq:l_step}}
\State $\widehat{\mathbf S}_{k+1} \gets \widehat{\mathbf U}_{k+1}^{\top} ({\mathbf Y}_k + h\,{\mathbf F}_k) \widehat{\mathbf V}_{k+1}$; \Comment{S-step \eqref{eq:s_step}}
\State $\widehat{\mathbf Y}_{k+1} \gets \widehat{\mathbf U}_{k+1} \widehat{\mathbf S}_{k+1} \widehat{\mathbf V}_{k+1}^\top$; 
\State ${\mathbf Y}_{k+1} \gets \typo{\Ro_{\M_r}}\bigl(\widehat{\mathbf Y}_{k+1}\bigr)$. \Comment{Truncation step \eqref{eq:t_step}}
\EndFor
\end{algorithmic}
\end{algorithm}

\subsection{\modif{Extension to high-order explicit Runge--Kutta methods}}
\label{sec:rk--bug}

\modif{Given the splitting-free reinterpretation, we now extend the BUG integrator to higher-order} explicit Runge--Kutta methods:
\begin{equation}
\label{eq:rk}
\begin{aligned}
\typo{{\mathbf X}}_{k i} &= \typo{{\mathbf X}}_{k} + h \sum_{j=1}^{i-1} a_{i j} {\mathbf F}(t_{k j},\typo{{\mathbf X}}_{k j}),
\quad i = 1, \ldots, s, \\
\typo{{\mathbf X}}_{k+1} &= \typo{{\mathbf X}}_{k} + h \sum_{i=1}^s b_{i} {\mathbf F}(t_{k i},\typo{{\mathbf X}}_{k i}),
\end{aligned}
\end{equation}
where $t_{k i} = t_k + c_i h$. The main idea is to \typo{apply one BUG step} at each stage of the Runge--Kutta method. \modif{The only point that requires special care is the construction of the augmented bases. For example, consider the projected discrete solution of~\eqref{eq:projected_ode} at the final stage starting from ${\mathbf Y}_k$:
\begin{equation*}
\begin{aligned}
\widetilde{\mathbf Z}_{k+1}
  =&\; {\mathbf Y}_k
   + h \sum_{i=1}^s b_{i}\,
   \Po_{{\mathbf Y}_{k i}}\!\bigl({\mathbf F}_{k i}\bigr) \\
  =&\; {\mathbf U}_k {\mathbf S}_k {\mathbf V}_k^{\top}
  + h \sum_{i=1}^s b_{i} (
       {\mathbf U}_{k i}{\mathbf U}_{k i}^{\top}{\mathbf F}_{k i}
     - {\mathbf U}_{k i}{\mathbf U}_{k i}^{\top}{\mathbf F}_{k i}{\mathbf V}_{k i}{\mathbf V}_{k i}^{\top}
     + {\mathbf F}_{k i}{\mathbf V}_{k i}{\mathbf V}_{k i}^{\top} ),
\end{aligned}
\end{equation*}
where ${\mathbf F}_{k i} = {\mathbf F}(t_{k i}, {\mathbf Y}_{k i})$. In order to provide an exact representation of the projected discrete solution $\widetilde{\mathbf Z}_{k+1}$, the augmented bases $(\widehat{\mathbf U}_{k+1}, \widehat{\mathbf V}_{k+1})$ are chosen such that they span the column and row spaces of ${\mathbf Y}_k$ and of the tangent-space projections $\Po_{{\mathbf Y}_{k i}}\!\bigl({\mathbf F}_{k i}\bigr)$ for $i \in \{1,\ldots,s\}$, that is,
\begin{equation*}
\begin{aligned}
  \widehat{\mathbf U}_{k+1}
  &= \operatorname{orth}(
      [\,{\mathbf U}_k,\;
      {\mathbf U}_{k 1},\; {\mathbf F}_{k 1}{\mathbf V}_{k 1},\;
      \ldots,\;
      {\mathbf U}_{k s},\; {\mathbf F}_{k s}{\mathbf V}_{k s}\,]
      ), \\
  \widehat{\mathbf V}_{k+1}
  &= \operatorname{orth}(
      [\,{\mathbf V}_k,\;
      {\mathbf V}_{k 1},\; {\mathbf F}_{k 1}^{\top}{\mathbf U}_{k 1},\;
      \ldots,\;
      {\mathbf V}_{k s},\; {\mathbf F}_{k s}^{\top}{\mathbf U}_{k s}\,]
      ).
\end{aligned}
\end{equation*}
Since ${\mathbf U}_{k 1} = {\mathbf U}_k$ and ${\mathbf V}_{k 1} = {\mathbf V}_k$, these terms can be removed from the concatenations above. Moreover, to exclude directions that do not contribute to the final combination, the coefficients $b_i$ are kept so that the corresponding bases are automatically discarded when $b_i = 0$. Hence, the augmented bases are finally defined as
\begin{equation*}
\begin{aligned}
  \widehat{\mathbf U}_{k+1}
  &= \operatorname{orth}(
      [\,{\mathbf U}_k,\;
      b_1 {\mathbf F}_{k 1}{\mathbf V}_{k 1},\;
      \ldots,\;
      b_s {\mathbf U}_{k s},\; b_s {\mathbf F}_{k s}{\mathbf V}_{k s}\,]
      ), \\
  \widehat{\mathbf V}_{k+1}
  &= \operatorname{orth}(
      [\,{\mathbf V}_k,\;
      b_1 {\mathbf F}_{k 1}^{\top}{\mathbf U}_{k 1},\;
      \ldots,\;
      b_s {\mathbf V}_{k s},\; b_s {\mathbf F}_{k s}^{\top}{\mathbf U}_{k s}\,]
      ).
\end{aligned}
\end{equation*}
The complete Runge--Kutta BUG (RK--BUG) integrator is summarized in Algorithm~\ref{alg:RKBUG}.} \review{Note that the rank $r$ is fixed here for simplicity, but an adaptive rank can also be used to truncate the augmented solution (see Section~\ref{sec:rank_adaptive} for a rank-adaptive strategy).}

\begin{algorithm}[H]
\caption{\typo{Runge--Kutta Basis-Update \& Galerkin (RK--BUG) integrator}}
\label{alg:RKBUG}
\begin{algorithmic}[1]
\Require {${\mathbf Y}_0$.}
\Ensure {${\mathbf Y}_{1},\ldots, {\mathbf Y}_{N}$.}
\For {$k=0,\ldots,N-1$}
\State ${\mathbf Y}_{k 1} \gets {\mathbf Y}_k$;
\For {$i=1,\ldots,s-1$}
\State ${\mathbf Y}_{k i} := {\mathbf U}_{k i} {\mathbf S}_{k i} {\mathbf V}_{k i}^\top$;
\State ${\mathbf F}_{k i} \gets {\mathbf F}(t_k + c_i h,{\mathbf Y}_{k i})$;
\State $\widehat{\mathbf U}_{k, i+1} \gets \operatorname{orth}([\,{\mathbf U}_k,\;\modif{a_{i+1,1}} {\mathbf F}_{k 1}{\mathbf V}_{k 1},\;\ldots,\;\modif{a_{i+1,i}} {\mathbf U}_{k i},\; \modif{a_{i+1,i}} {\mathbf F}_{k i}{\mathbf V}_{k i}\,])$; 
\State $\widehat{\mathbf V}_{k,i+1} \gets \operatorname{orth}([\,{\mathbf V}_k,\;\modif{a_{i+1,1}} {\mathbf F}_{k 1}^{\top}{\mathbf U}_{k 1},\;\ldots,\;\modif{a_{i+1,i}} {\mathbf V}_{k i},\; \modif{a_{i+1,i}} {\mathbf F}_{k i}^{\top}{\mathbf U}_{k i}\,])$; 
\State $\widehat{\mathbf S}_{k,i+1} \gets \widehat{\mathbf U}_{k,i+1}^\top \bigl( {\mathbf Y}_{k} + h \, ( a_{i+1,1} {\mathbf F}_{k 1} + \ldots + a_{i+1,i} {\mathbf F}_{k i} ) \bigr) \widehat{\mathbf V}_{k,i+1}$; 
\State $\widehat{\mathbf Y}_{k,i+1} \gets \widehat{\mathbf U}_{k,i+1} \widehat{\mathbf S}_{k,i+1} \widehat{\mathbf V}_{k,i+1}^\top$;
\State ${\mathbf Y}_{k,i+1} \gets \typo{\Ro_{\M_r}}\bigl(\widehat{\mathbf Y}_{k,i+1} \bigr)$; 
\EndFor
\State ${\mathbf Y}_{k s} := {\mathbf U}_{k s} {\mathbf S}_{k s} {\mathbf V}_{k s}^\top$;
\State ${\mathbf F}_{k s} \gets {\mathbf F}(t_k + c_s h,{\mathbf Y}_{k s})$;
\State $\widehat{\mathbf U}_{k+1} \gets \operatorname{orth}([\,{\mathbf U}_k,\;\modif{b_{1}} {\mathbf F}_{k 1}{\mathbf V}_{k 1},\;\ldots,\;\modif{b_{s}} {\mathbf U}_{k s},\; \modif{b_{s}} {\mathbf F}_{k s}{\mathbf V}_{k s}\,])$; 
\State $\widehat{\mathbf V}_{k+1} \gets \operatorname{orth}([\,{\mathbf V}_k,\;\modif{b_1} {\mathbf F}_{k 1}^{\top}{\mathbf U}_{k 1},\;\ldots,\;\modif{b_s} {\mathbf V}_{k s},\; \modif{b_s} {\mathbf F}_{k s}^{\top}{\mathbf U}_{k s}\,])$; 
\State $\widehat{\mathbf S}_{k+1} \gets \widehat{\mathbf U}_{k+1}^\top \bigl( {\mathbf Y}_{k} + h \, ( b_{1} {\mathbf F}_{k 1} + \ldots + b_{s} {\mathbf F}_{k s} ) \bigr) \widehat{\mathbf V}_{k+1}$; 
\State $\widehat{\mathbf Y}_{k+1} \gets \widehat{\mathbf U}_{k+1} \widehat{\mathbf S}_{k+1} \widehat{\mathbf V}_{k+1}^\top$; 
\State ${\mathbf Y}_{k+1} \gets \typo{\Ro_{\M_r}}\bigl(\widehat{\mathbf Y}_{k+1} \bigr)$. 
\EndFor
\end{algorithmic}
\end{algorithm}

\subsection{\review{Conservative variant}}

\review{
The matrix differential equation~\eqref{eq:ode} may admit local conservation laws. Suppose for example that the locally conserved quantities are obtained by right-projecting the solution onto fixed directions ${\mathbf W} \in \R^{m \times r_{\rm cons}}$, that is ${\mathbf M}(t) = {\mathbf X}(t)\,{\mathbf W} \in \R^{n \times r_{\rm cons}}$, and that the corresponding global invariants are recovered by left-multiplying ${\mathbf M}(t)$ by a constant vector ${\mathbf w}_x \in \R^n$ (usually representing the spatial averaging), so that $\langle {\mathbf M}(t)\rangle_x = {\mathbf w}_x^\top {\mathbf M}(t)$. Since ${\mathbf X}(t)$ is the solution of~\eqref{eq:ode}, the projected quantities ${\mathbf M}(t)$ satisfy
\begin{equation}
\label{eq:local_conservation}
\dot{\mathbf M}(t) = {\mathbf F}(t,{\mathbf X}(t))\,{\mathbf W}.
\end{equation}
In kinetic equations such as Boltzmann--BGK or Vlasov--Poisson, the solution ${\mathbf X}(t)$ typically represents the distribution function $f({\boldsymbol x},{\boldsymbol v},t)$, with spatial nodes stored along the rows and velocity nodes along the columns. The corresponding locally conserved quantities are the mass, momentum, and energy per unit volume:
\begin{equation}
\label{eq:moments}
[\rho({\boldsymbol x},t),\, {\boldsymbol j}({\boldsymbol x},t),\, E({\boldsymbol x},t)] = \int_{\Omega_v} f({\boldsymbol x},{\boldsymbol v},t)\,\bigl[1,\, {\boldsymbol v},\, \tfrac12\|{\boldsymbol v}\|_2^2\bigr]\,\diff{\boldsymbol v}.
\end{equation}
At the discrete level, the analogue of \eqref{eq:moments} can be written as ${\mathbf M}(t) = {\mathbf X}(t)\,{\mathbf W}$, where ${\mathbf W}$ incorporates both the evaluations of $[1,\,\boldsymbol v,\,\tfrac12\|\boldsymbol v\|_2^2]$ on the velocity grid and the quadrature weights used in the numerical integration. Under periodic boundary conditions in ${\boldsymbol x}$ or vanishing fluxes at the boundary, the spatial averages $\langle {\mathbf M}(t)\rangle_x$ become invariants of the dynamics, since they satisfy $\frac{\diff}{\diff t}\,\langle {\mathbf M}(t)\rangle_x = {\mathbf 0}$ when the spatial and velocity discretizations preserve the conservation properties.
}

\review{
The RK--BUG integrator can be modified at little additional computational cost to maintain such conservation properties. The idea is to decompose the solution into a conservative component and a remainder:
\begin{equation*}
{\mathbf Y}(t) = {\mathbf K}(t)\,{\mathbf V}_{\mathrm{cons}}^\top + {\mathbf U}(t)\,{\mathbf S}(t)\,{\mathbf V}(t)^\top,
\end{equation*}
where ${\mathbf V}_{\rm cons} := \operatorname{ortho}({\mathbf W})$ contains the conservative modes, and ${\mathbf V}(t)$ is constrained to remain orthogonal to ${\mathbf V}_{\rm cons}$ for all $t$. The evolution of the low-rank factors $({\mathbf K},{\mathbf U},{\mathbf S},{\mathbf V})$ is governed by
\begin{subequations}
\label{eq:cons_evo}
\begin{empheq}[left={\empheqlbrace\,}]{align}
\dot{\mathbf K} &= {\mathbf F}(t,{\mathbf Y})\,{\mathbf V}_{\mathrm{cons}}, \label{eq:K_evo}\\
\frac{\diff}{\diff t}\bigl({\mathbf U}{\mathbf S}{\mathbf V}^\top\bigr) &= \Po_{{\mathbf U}{\mathbf S}{\mathbf V}^\top}\bigl({\mathbf F}(t,{\mathbf Y}) - {\mathbf F}(t,{\mathbf Y}) {\mathbf V}_{\mathrm{cons}}{\mathbf V}_{\mathrm{cons}}^\top\bigr), \label{eq:rem_evo}
\end{empheq}
\end{subequations}
where the orthonormality constraint ${\mathbf V}_{\mathrm{cons}}^\top \dot{\mathbf V} = {\mathbf 0}$ is enforced in~\eqref{eq:rem_evo} by subtracting the term ${\mathbf F}(t,{\mathbf Y})\,{\mathbf V}_{\mathrm{cons}}{\mathbf V}_{\mathrm{cons}}^\top$, which ensures that ${\mathbf V}$ remains orthogonal to ${\mathbf V}_{\mathrm{cons}}$ for all times. For the initialization, let ${\mathbf X}_0^\perp := {\mathbf X}_0\bigl({\mathbf I} - {\mathbf V}_{\mathrm{cons}}{\mathbf V}_{\mathrm{cons}}^\top\bigr)$ denote the projection of the initial condition onto the orthogonal complement of ${\mathbf V}_{\rm cons}$. Given the truncated SVD ${\mathbf X}_0^\perp \approx {\mathbf \Phi}{\mathbf \Sigma}{\mathbf \Psi}^\top$, the initial factors are
\begin{equation*}
{\mathbf K}_0 = {\mathbf X}_0{\mathbf V}_{\mathrm{cons}}, \qquad
{\mathbf U}_0 = {\mathbf \Phi}, \qquad
{\mathbf S}_0 = {\mathbf \Sigma}, \qquad
{\mathbf V}_0 = {\mathbf \Psi}.
\end{equation*}
Then, system~\eqref{eq:cons_evo} is integrated using an explicit Runge--Kutta method: the conservative equation~\eqref{eq:K_evo} is advanced by the underlying Runge--Kutta scheme, while the remainder equation~\eqref{eq:rem_evo} is updated via the corresponding RK--BUG integrator. In practice, although ${\mathbf V}(t)$ should theoretically remain orthogonal to ${\mathbf V}_{\rm cons}$, an additional orthonormalization step is applied to enforce that the augmented basis $\widehat{\mathbf V}$ satisfies ${\mathbf V}_{\rm cons}^\top \widehat{\mathbf V} = 0$ and, consequently, the orthonormality constraint ${\mathbf V}_{\rm cons}^\top {\mathbf V} = 0$ holds up to machine precision. Moreover, the RK--BUG integrator is applied directly to ${\mathbf F}(t,{\mathbf Y})$, rather than to ${\mathbf F}(t,{\mathbf Y})\bigl({\mathbf I} - {\mathbf V}_{\mathrm{cons}}{\mathbf V}_{\mathrm{cons}}^\top\bigr)$, since the contribution ${\mathbf F}(t,{\mathbf Y}) {\mathbf V}_{\mathrm{cons}}{\mathbf V}_{\mathrm{cons}}^\top$ vanishes when projected onto ${\mathbf V}$ or during the orthonormalization step. The full procedure is summarized in Algorithm~\ref{alg:consRKBUG}. With this construction, the conservative component captures the dynamics of the local conservation law~\eqref{eq:local_conservation}, ensuring that the corresponding conserved quantities (and thus the associated global invariants) are preserved.
}

\begin{algorithm}[H]
\caption{\review{Conservative RK--BUG integrator}}
\label{alg:consRKBUG}
\review{
\begin{algorithmic}[1]
\Require {${\mathbf Y}_0$.}
\Ensure {${\mathbf Y}_{1},\ldots, {\mathbf Y}_{N}$.}
\For {$k=0,\ldots,N-1$}
\State ${\mathbf Y}_{k 1} \gets {\mathbf Y}_k$;
\For {$i=1,\ldots,s-1$}
\State ${\mathbf Y}_{k i} := {\mathbf K}_{ki} {\mathbf V}_{\mathrm{cons}}^\top + {\mathbf U}_{k i} {\mathbf S}_{k i} {\mathbf V}_{k i}^\top$;
\State ${\mathbf F}_{k i} \gets {\mathbf F}(t_k + c_i h,{\mathbf Y}_{k i})$;
\State ${\mathbf K}_{k, i+1} \gets {\mathbf K}_{k} + h \, ( a_{i+1,1} {\mathbf F}_{k 1} + \ldots + a_{i+1,i} {\mathbf F}_{k i} ) {\mathbf V}_{\mathrm{cons}}$;
\State $\widehat{\mathbf U}_{k, i+1} \gets \operatorname{orth}([\,{\mathbf U}_k,\;a_{i+1,1} {\mathbf F}_{k 1} {\mathbf V}_{k 1},\;\ldots,\;a_{i+1,i} {\mathbf U}_{k i},\; a_{i+1,i} {\mathbf F}_{k i} {\mathbf V}_{k i}\,])$; 
\State $\widehat{\mathbf V}_{k,i+1} \gets \operatorname{orth}([\,{\mathbf V}_k,\;a_{i+1,1} {\mathbf F}_{k 1}^{\top}{\mathbf U}_{k 1},\;\ldots,\;a_{i+1,i} {\mathbf V}_{k i},\; a_{i+1,i} {\mathbf F}_{k i}^{\top}{\mathbf U}_{k i}\,])$; 
\State $[\,\sim,\;\widehat{\mathbf V}_{k,i+1}\,] \gets \operatorname{orth}([\,{\mathbf V}_{\mathrm{cons}},\;\widehat{\mathbf V}_{k,i+1}\,])$; 
\State $\widehat{\mathbf S}_{k,i+1} \gets \widehat{\mathbf U}_{k,i+1}^\top \bigl( {\mathbf U}_{k} {\mathbf S}_{k} {\mathbf V}_{k}^\top + h \, ( a_{i+1,1} {\mathbf F}_{k 1} + \ldots + a_{i+1,i} {\mathbf F}_{k i} ) \bigr) \widehat{\mathbf V}_{k,i+1}$; 
\State $\widehat{\mathbf Y}_{k,i+1}^{\mathrm{rem}} \gets \widehat{\mathbf U}_{k,i+1} \widehat{\mathbf S}_{k,i+1} \widehat{\mathbf V}_{k,i+1}^\top$;
\State ${\mathbf Y}_{k,i+1} \gets {\mathbf K}_{k,i+1} {\mathbf V}_{\mathrm{cons}}^\top + \Ro_{\M_r}\bigl( \widehat{\mathbf Y}_{k,i+1}^{\mathrm{rem}} \bigr)$; 
\EndFor
\State ${\mathbf Y}_{k s} := {\mathbf K}_{k s} {\mathbf V}_{\mathrm{cons}}^\top + {\mathbf U}_{k s} {\mathbf S}_{k s} {\mathbf V}_{k s}^\top$;
\State ${\mathbf F}_{k s} \gets {\mathbf F}(t_k + c_s h,{\mathbf Y}_{k s})$;
\State ${\mathbf K}_{k+1} \gets {\mathbf K}_{k} + h \, ( b_{1} {\mathbf F}_{k 1} + \ldots + b_{s} {\mathbf F}_{k s} ) {\mathbf V}_{\mathrm{cons}}$;
\State $\widehat{\mathbf U}_{k+1} \gets \operatorname{orth}([\,{\mathbf U}_k,\;b_{1} {\mathbf F}_{k 1} {\mathbf V}_{k 1},\;\ldots,\;b_{s} {\mathbf U}_{k s},\; b_{s} {\mathbf F}_{k s} {\mathbf V}_{k s}\,])$; 
\State $\widehat{\mathbf V}_{k+1} \gets \operatorname{orth}([\,{\mathbf V}_k,\;b_1 {\mathbf F}_{k 1}^{\top}{\mathbf U}_{k 1},\;\ldots,\;b_s {\mathbf V}_{k s},\; b_s {\mathbf F}_{k s}^{\top}{\mathbf U}_{k s}\,])$;
\State $[\,\sim,\;\widehat{\mathbf V}_{k+1}\,] \gets \operatorname{orth}([\,{\mathbf V}_{\mathrm{cons}},\;\widehat{\mathbf V}_{k+1}\,])$;
\State $\widehat{\mathbf S}_{k+1} \gets \widehat{\mathbf U}_{k+1}^\top \bigl( {\mathbf U}_{k} {\mathbf S}_{k} {\mathbf V}_{k}^\top + h \, ( b_{1} {\mathbf F}_{k 1} + \ldots + b_{s} {\mathbf F}_{k s} ) \bigr) \widehat{\mathbf V}_{k+1}$; 
\State $\widehat{\mathbf Y}_{k+1}^{\mathrm{rem}} \gets \widehat{\mathbf U}_{k+1} \widehat{\mathbf S}_{k+1} \widehat{\mathbf V}_{k+1}^\top$; 
\State ${\mathbf Y}_{k+1} \gets {\mathbf K}_{k+1} {\mathbf V}_{\mathrm{cons}}^\top + \Ro_{\M_r}\bigl( \widehat{\mathbf Y}_{k+1}^{\mathrm{rem}} \bigr)$. 
\EndFor
\end{algorithmic}
}
\end{algorithm}

\subsection{\review{Computational cost}}

\review{We now analyze the computational cost of the RK--BUG integrator. For simplicity, we exclude the cost of evaluating~$\mathbf F$, more precisely its 
projection onto the appropriate bases, since it depends on the particular form of~$\mathbf F$ (e.g., linear, polynomial, or nonlinear) and not on the time-integration method itself.}

\review{The memory footprint is $\mathcal O(s\,(n{+}m)\,r)$, corresponding to the storage of the bases $\{({\mathbf U}_{ki},{\mathbf V}_{ki})\}_{i=1}^s$, with a temporary buffer of size $\mathcal O((n{+}m)\,\hat r + \hat r^{2})$ for the current augmented factors $(\widehat{\mathbf U}_{ki},\widehat{\mathbf S}_{ki},\widehat{\mathbf V}_{ki})$.}

\review{The complexity of Algorithm~\ref{alg:RKBUG} is dominated at each step by three operations: (1) \emph{thin QR factorizations} of $n\times \hat r$ and $m\times \hat r$ matrices, with cost $\mathcal O\!\big(s\,(n{+}m)\,\hat r^{2}\big)$, (2) \emph{matrix--matrix multiplications} involving reduced factors, with cost $\mathcal O\!\big(s\,(n{+}m)\,r\,\hat r\big)$, and (3) \emph{truncated SVD} on $\hat r\times \hat r$ blocks, with cost $\mathcal O\!\big(s\,\hat r^{3}\big)$. Hence, the overall complexity per step is}
\begin{equation*}
\review{\mathcal O\!\big(s\,(n{+}m)\,\hat r^{2}\big).}
\end{equation*}

\review{It is important to note that the computational efficiency of the RK--BUG integrator depends on the augmented rank $\hat r$. To mitigate its growth, two key components are employed in Algorithm~\ref{alg:RKBUG}. First, the coefficients $(a_{ij},b_i)$ are kept explicitly in the construction of the augmented bases, so that the corresponding bases are discarded when these coefficients are zero. Second, a truncation step is applied after each stage. As a result, the augmented rank $\hat r$ is at most $2ir$ at stage~$i$ and $2sr$ at the final stage.}

\review{
The conservative RK--BUG integrator (Algorithm~\ref{alg:consRKBUG}) introduces 
only minor additional operations compared with its non-conservative variant. The extra cost comes from: (1) \emph{the update of ${\mathbf K}$}, with complexity $\mathcal{O}(s\,m\,r_{\rm cons})$, and (2) \emph{the additional orthonormalization step}, with complexity $\mathcal{O}\big(s\,m\,(r_{\rm cons}+ \hat r)^2\big)$. As a result, the overall per-step complexity becomes
\begin{equation*}
\mathcal{O}\!\bigl(s\,n\,\hat r^{2} \,+\, s\,m\,(r_{\rm cons}+\hat r)^{2}\bigr),
\end{equation*}
which remains of the same order of magnitude as in the non-conservative case.}

\section{Convergence analysis}
\label{sec:3}

In this section, we \typo{establish} the \typo{high-order} convergence of the \typo{RK--BUG} integrator. \typo{The analysis is conducted} for a step size $h \leq h_0$ \typo{(small enough, see~\cite{harrier1993solving})}, and for a fixed rank~$r$, \review{in the non-conservative case}. \typo{In the following, all estimates are expressed} in the Frobenius norm~$\norm{\cdot}_F$, and \typo{we denote by}~$\dott{\cdot}{\cdot}_F$ the \typo{associated} inner product.

\subsection{\typo{Assumptions and} preliminary results}

The \typo{convergence} analysis \typo{relies} on two assumptions. The first \typo{guarantees} the existence and uniqueness of the exact solution~\typo{${\mathbf X}$} according to the Picard--Lindel\"of theorem, while the second \typo{is a standard assumption ensuring the high-order convergence of explicit} Runge--Kutta methods \typo{(see Theorem~3.1 in Chapter~2 of~\cite{harrier1993solving})}. Compared to \typo{the convergence analysis of the BUG integrator}~\cite{ceruti2022rank}, we do not assume that ${\mathbf F}(t,{\mathbf Z})$ is bounded for all ${\mathbf Z} \in \R^{n \times m}$. \modif{Instead, starting from Assumption~\ref{theo:A1}, we prove in Propositions~\ref{theo:P1} and~\ref{theo:P2} that the solution remains in a neighbourhood of the initial condition ${\mathbf Y}_0$. It follows from the Lipschitz continuity of ${\mathbf F}$ that ${\mathbf F}$ is bounded on this compact set, which is sufficient for the subsequent analysis.} \typo{Finally}, we show \typo{in Lemma~\ref{theo:L4}} that the exact flow ${\boldsymbol\Phi}^{t}_{\mathbf F}$ is Lipschitz continuous, \typo{as} a direct consequence of \typo{Assumption~\ref{theo:A1}}.

\medskip

\begin{assumption}
\label{theo:A1}
${\mathbf F}(t,{\mathbf Z})$ is continuous in time and Lipschitz continuous in ${\mathbf Z}$: there exists a Lipschitz constant $L>0$ (independent of $t$) such that
\begin{equation*}
\norm{{\mathbf F}(t,{\mathbf Z}_1)-{\mathbf F}(t,{\mathbf Z}_2)}_F \leq L \norm{{\mathbf Z}_1-{\mathbf Z}_2}_F,
\end{equation*}
for all $t \in \ff{0}{T}$ and ${\mathbf Z}_1,{\mathbf Z}_2 \in \R^{n \times m}$.
\end{assumption}

\medskip

\begin{assumption}
\label{theo:A2}
\typo{Let $(a_{ij},b_j,c_i)$ define an explicit Runge--Kutta method~\eqref{eq:rk} of order~$p$.} The first $p$ derivatives of the exact flow ${\boldsymbol\Phi}^t_{\mathbf F}(\typo{{\mathbf Z}})$ (i.e., the mapping such that $\typo{{\mathbf X}}(t) = {\boldsymbol\Phi}^t_{\mathbf F}(\typo{{\mathbf X}}_0)$, \typo{with ${\mathbf X}(t)$ the exact solution of~\eqref{eq:ode}}) exist and are continuous for all $\typo{{\mathbf Z}} \in \R^{n \times m}$.
\end{assumption}

\medskip

\begin{proposition}
\label{theo:P1}
Suppose that Assumption~\ref{theo:A1} holds. \typo{Then}, the exact solution ${\mathbf Y}(t)$ of the projected differential equation~\eqref{eq:projected_ode} \typo{satisfies}
\begin{equation}
\norm{{\mathbf Y}(t)-{\mathbf Y}_0}_F \leq  \int_0^t e^{L(t-s)} \norm{{\mathbf F}(s,{\mathbf Y}_0)}_F \, \mathrm{d}s,
\end{equation}
for all $t \in \ff{0}{T}$.
\end{proposition}
\begin{proof}
Without loss of generality, assume that ${\mathbf Y}(t) \neq {\mathbf Y}_0$ for all $t \in \oo{0}{T}$. If ${\mathbf Y}(t) = {\mathbf Y}_0$ for certain times, \typo{the analysis can be carried out} independently on the \typo{different} subintervals where ${\mathbf Y}(t) \neq {\mathbf Y}_0$. \typo{From Assumption~\ref{theo:A1}, we obtain} the differential inequality
\begin{equation*}
\begin{aligned}
\norm{{\mathbf Y}(t)-{\mathbf Y}_0}_F \frac{\diff}{\diff t} \norm{{\mathbf Y}(t)-{\mathbf Y}_0}_F &= \frac 1 2 \frac{\diff}{\diff t} \norm{{\mathbf Y}(t)-{\mathbf Y}_0}_F^2 \\
&= \dott{{\mathbf Y}(t)-{\mathbf Y}_0}{\Po_{\mathbf Y}\big({\mathbf F}(t,{\mathbf Y}(t))\big)}_F \\
&\leq \norm{{\mathbf Y}(t)-{\mathbf Y}_0}_F \, \norm{\Po_{\mathbf Y}\big({\mathbf F}(t,{\mathbf Y}(t))\big)}_F \\
&\leq \norm{{\mathbf Y}(t)-{\mathbf Y}_0}_F \, \norm{{\mathbf F}(t,{\mathbf Y}(t))}_F \\
&\leq \norm{{\mathbf Y}(t)-{\mathbf Y}_0}_F \big(\norm{{\mathbf F}(t,{\mathbf Y}(t))-{\mathbf F}(t,{\mathbf Y}_0)}_F + \norm{{\mathbf F}(t,{\mathbf Y}_0)}_F \big) \\
&\leq \norm{{\mathbf Y}(t)-{\mathbf Y}_0}_F \big( L \norm{{\mathbf Y}(t)-{\mathbf Y}_0}_F + \norm{{\mathbf F}(t,{\mathbf Y}_0)}_F \big),
\end{aligned}
\end{equation*}
which can be rewritten as
\begin{equation*}
\frac{\diff}{\diff t} \norm{{\mathbf Y}(t)-{\mathbf Y}_0}_F \leq L \norm{{\mathbf Y}(t)-{\mathbf Y}_0}_F + \norm{{\mathbf F}(t,{\mathbf Y}_0)}_F.
\end{equation*}
\typo{Then}, according to Gr\"onwall's inequality, the exact solution ${\mathbf Y}(t)$ \typo{satisfies}
\begin{equation*}
\norm{{\mathbf Y}(t)-{\mathbf Y}_0}_F \leq \int_0^t e^{L(t-s)} \norm{{\mathbf F}(s,{\mathbf Y}_0)}_F \, \mathrm{d}s,
\end{equation*}
which concludes the proof.
\end{proof}

\medskip

\begin{lemma}
\label{theo:L1}
Suppose that Assumption~\ref{theo:A1} holds. Then, for $h \leq h_0$ \modif{and a fixed rank~$r$}, the \typo{RK--BUG} solution ${\mathbf Y}_{k i}$ satisfies
\begin{equation}
\norm{{\mathbf Y}_{k i}-{\mathbf Y}_0}_F 
   \leq \big( 1+K_{i 0} h \big) \norm{{\mathbf Y}_{k}-{\mathbf Y}_0}_F 
      + h \sum_{j=1}^{i-1} K_{i j}\norm{{\mathbf F}(t_{k j},{\mathbf Y}_0)}_F,
\end{equation}
\typo{at each stage $i \in \{1,\ldots,s\}$. Here}, the constants $K_{i j} \geq 0$ are independent of $h$ \modif{and $r$}.
\end{lemma}
\begin{proof}
We proceed by induction \typo{on the stage index $i$}. For $i=1$, the statement is trivial with \modif{$K_{1,0}=0$}, since ${\mathbf Y}_{k 1} = {\mathbf Y}_{k}$. \typo{Assume now that the result holds for all stages $j < i$, with $i \in \{2,\ldots,s\}$}. Then, \typo{we have}
\begin{equation*}
\begin{aligned}
\norm{{\mathbf Y}_{k i}-{\mathbf Y}_0}_F 
&= \norm{ \typo{\Ro_{\M_r}}\bigl(\widehat{\mathbf Y}_{k i} \bigr) - {\mathbf Y}_0}_F \\
&\leq \norm{\typo{\Ro_{\M_r}}\bigl( \widehat{\mathbf Y}_{k i} \bigr) - \widehat{\mathbf Y}_{k i}}_F + \norm{\widehat{\mathbf Y}_{k i} - {\mathbf Y}_0}_F \\
&= \min_{{\mathbf Z} \in \M_r} \norm{{\mathbf Z} - \widehat{\mathbf Y}_{k i}}_F + \norm{\widehat{\mathbf Y}_{k i} - {\mathbf Y}_0}_F \\
&\leq \norm{{\mathbf Y}_{k} - \widehat{\mathbf Y}_{k i}}_F + \norm{\widehat{\mathbf Y}_{k i} - {\mathbf Y}_0}_F \\
&\leq \norm{{\mathbf Y}_{k}-{\mathbf Y}_0}_F + 2h \sum_{j=1}^{i-1} |a_{i j}|\, \norm{\widehat{\mathbf U}_{k i}\widehat{\mathbf U}_{k i}^\top {\mathbf F}(t_{k j},{\mathbf Y}_{k j}) \widehat{\mathbf V}_{k i}\widehat{\mathbf V}_{k i}^\top}_F \\
&\leq \norm{{\mathbf Y}_{k}-{\mathbf Y}_0}_F + 2h \sum_{j=1}^{i-1} |a_{i j}|\,\norm{{\mathbf F}(t_{k j},{\mathbf Y}_{k j})}_F \\
&\leq \norm{{\mathbf Y}_{k}-{\mathbf Y}_0}_F + 2h \sum_{j=1}^{i-1} |a_{i j}|\, \Big( \norm{{\mathbf F}(t_{k j},{\mathbf Y}_{k j})-{\mathbf F}(t_{k j},{\mathbf Y}_0)}_F + \norm{{\mathbf F}(t_{k j},{\mathbf Y}_0)}_F \Big) \\
&\leq \norm{{\mathbf Y}_{k}-{\mathbf Y}_0}_F + 2h \sum_{j=1}^{i-1} |a_{i j}|\, \Big( L \norm{{\mathbf Y}_{k j}-{\mathbf Y}_0}_F + \norm{{\mathbf F}(t_{k j},{\mathbf Y}_0)}_F \Big) \\
&\leq \Big( 1 + 2h \sum_{j=1}^{i-1} |a_{i j}| L (1+K_{j 0}h ) \Big) 
        \norm{{\mathbf Y}_{k}-{\mathbf Y}_0}_F \\
&\quad + 2h \sum_{j=1}^{i-1} |a_{i j}|\, 
        \Big( \norm{{\mathbf F}(t_{k j},{\mathbf Y}_0)}_F 
             + Lh \sum_{l=1}^{j-1} K_{j l} 
               \norm{{\mathbf F}(t_{k l},{\mathbf Y}_0)}_F \Big) \\
&\leq \big( 1 + K_{i 0} h \big) \norm{{\mathbf Y}_{k}-{\mathbf Y}_0}_F 
      + h \sum_{l=1}^{i-1} K_{i l} \norm{{\mathbf F}(t_{k l},{\mathbf Y}_0)}_F,
\end{aligned}
\end{equation*}
where
\begin{equation*}
\begin{aligned}
K_{i 0} &= 2L \sum_{j=1}^{i-1} |a_{i j}| (1+K_{j 0}h_0), \\
K_{i l} &= 2 |a_{i l}| + 2 L h_0 \sum_{j=l+1}^{i-1} |a_{i j}| K_{j l}, 
           \quad l=1,\ldots,i-2, \\
K_{i,i-1} &= 2 |a_{i,i-1}|.
\end{aligned}
\end{equation*}
\end{proof}

\medskip

\begin{lemma}
\label{theo:L2}
Suppose that Assumption~\ref{theo:A1} holds. Then, for $h \leq h_0$ \modif{and a fixed rank~$r$}, the \typo{RK--BUG} solution ${\mathbf Y}_{k+1}$ \typo{satisfies}
\begin{equation}
\norm{{\mathbf Y}_{k+1}-{\mathbf Y}_0}_F \leq \big( 1 + \widetilde K_{0} h \big) \norm{{\mathbf Y}_{k}-{\mathbf Y}_0}_F + h \sum_{i=1}^{s} \widetilde K_{i}\, \norm{{\mathbf F}(t_{k i},{\mathbf Y}_0)}_F,
\end{equation}
where the constants $\widetilde K_{i} \geq 0$ are independent of $h$ \modif{and $r$}.
\end{lemma}
\begin{proof}
\typo{We follow the same reasoning as in Lemma~\ref{theo:L1}}. \modif{For the initialization step ($k=0$), the statement is trivial since ${\mathbf Y}_{k} = {\mathbf Y}_0$. Assume now that the result holds for some $k \ge 0$}. Then, \typo{using} Lemma~\ref{theo:L1}, \typo{we obtain}
\begin{equation*}
\begin{aligned}
\norm{{\mathbf Y}_{k+1}-{\mathbf Y}_0}_F &\leq \norm{{\mathbf Y}_k - \widehat{\mathbf Y}_{k+1}}_F + \norm{\widehat{\mathbf Y}_{k+1} - {\mathbf Y}_0}_F \\
&\leq \norm{{\mathbf Y}_{k}-{\mathbf Y}_0}_F + 2h \sum_{i=1}^{s} |b_{i}|\, \norm{\widehat{\mathbf U}_{k s}\widehat{\mathbf U}_{k s}^\top {\mathbf F}(t_{k i},{\mathbf Y}_{k i}) \widehat{\mathbf V}_{k s}\widehat{\mathbf V}_{k s}^\top}_F \\
&\leq \norm{{\mathbf Y}_{k}-{\mathbf Y}_0}_F + 2h \sum_{i=1}^{s} |b_{i}|\, \norm{{\mathbf F}(t_{k i},{\mathbf Y}_{k i})}_F \\
&\leq \norm{{\mathbf Y}_{k}-{\mathbf Y}_0}_F + 2h \sum_{i=1}^{s} |b_{i}|\, \Big( \norm{{\mathbf F}(t_{k i},{\mathbf Y}_{k i})-{\mathbf F}(t_{k i},{\mathbf X}_0)}_F + \norm{{\mathbf F}(t_{k i},{\mathbf Y}_0)}_F \Big) \\
&\leq \norm{{\mathbf Y}_{k}-{\mathbf Y}_0}_F + 2h \sum_{i=1}^{s} |b_{i}|\, \Big( L \norm{{\mathbf Y}_{k i}-{\mathbf Y}_0}_F + \norm{{\mathbf F}(t_{k i},{\mathbf Y}_0)}_F \Big) \\
&\leq \Big( 1 + 2h \sum_{i=1}^{s} |b_{i}| L (1 + K_{i 0} h) \Big) \norm{{\mathbf Y}_{k}-{\mathbf Y}_0}_F \\
&\quad + 2h \sum_{i=1}^{s} |b_{i}|\, \Big( \norm{{\mathbf F}(t_{k i},{\mathbf Y}_0)}_F + Lh \sum_{j=1}^{i-1} K_{i j} \norm{{\mathbf F}(t_{k j},{\mathbf Y}_0)}_F \Big) \\
&\leq \big( 1 + \widetilde K_{0} h \big) \norm{{\mathbf Y}_{k}-{\mathbf Y}_0}_F + h \sum_{l=1}^{s} \widetilde K_{l}\, \norm{{\mathbf F}(t_{k l},{\mathbf Y}_0)}_F,
\end{aligned}
\end{equation*}
where 
\begin{equation*}
\begin{aligned}
\widetilde K_{0} &= 2L \sum_{i=1}^{s} |b_{i}| (1 + K_{i 0} h_0), \\
\widetilde K_{l} &= 2 |b_{l}| + 2 L h_0 \sum_{i=l+1}^{s} |b_{i}| K_{i l}, \quad l=1,\ldots,s-1, \\
\widetilde K_{s} &= 2 |b_{s}|.
\end{aligned}
\end{equation*}
\end{proof}

\medskip

\begin{proposition}
\label{theo:P2}
Suppose that Assumption~\ref{theo:A1} holds. Then, for $h \leq h_0$ \modif{and a fixed rank~$r$}, the \typo{RK--BUG} solution ${\mathbf Y}_{k}$ \typo{satisfies}
\begin{equation}
\norm{{\mathbf Y}_{k}-{\mathbf Y}_0}_F 
\leq  \frac{ e^{k h \widetilde K_{0}} - 1}{\widetilde K_{0}}
   \typo{\max_{0 \leq l \leq k-1}
   \Bigg( \sum_{i=1}^{s} \widetilde K_{i} 
   \norm{{\mathbf F}(t_{l i},{\mathbf Y}_0)}_F \Bigg)}.
\end{equation}
\end{proposition}
\begin{proof}
We \typo{first show} by induction that
\begin{equation}
\label{eq:P2_1}
\norm{{\mathbf Y}_{k}-{\mathbf Y}_0}_F \leq  h \sum_{l=0}^{k-1} \typo{G_l} \, (1 + h \widetilde K_{0})^{k-1-l},
\end{equation}
\typo{where $G_l = \sum_{i=1}^{s} \widetilde K_{i} \norm{{\mathbf F}(t_{l i},{\mathbf Y}_0)}_F$}. For $k=0$, the statement is trivial \typo{since ${\mathbf Y}_{k} = {\mathbf Y}_0$}. \typo{Assume now that the result holds for some $k \ge 0$}. Then, \typo{from} Lemma~\ref{theo:L2}, \typo{we obtain}
\begin{equation*}
\begin{aligned}
\norm{{\mathbf Y}_{k+1}-{\mathbf Y}_0}_F 
&\leq (1 + h \widetilde K_{0}) \norm{{\mathbf Y}_{k}-{\mathbf Y}_0}_F + h \typo{G_k} \\
&\leq (1 + h \widetilde K_{0}) \Big( h \sum_{l=0}^{k-1} \typo{G_l} \, (1 + h \widetilde K_{0})^{k-1-l}\Big) + h \typo{G_k} \\
&=  h \sum_{l=0}^{k} \typo{G_l} \, (1 + h \widetilde K_{0})^{k-l}.
\end{aligned}
\end{equation*}
Finally, \typo{using the induction result}~\eqref{eq:P2_1}, \typo{the solution ${\mathbf Y}_{k}$ satisfies}
\begin{equation*}
\begin{aligned}
\norm{{\mathbf Y}_{k}-{\mathbf Y}_0}_F &\leq \typo{\Big( \max_{0 \leq l \leq k-1} G_l \Big)} h \sum_{l=0}^{k-1} (1 + h \widetilde K_{0})^{k-1-l} \\
&=  \typo{\Big( \max_{0 \leq l \leq k-1} G_l \Big)} \frac{(1 + h \widetilde K_{0})^{k} - 1}{\widetilde K_{0}} \\
&\leq  \typo{\Big( \max_{0 \leq l \leq k-1} G_l \Big)} \frac{ e^{k h \widetilde K_{0}} - 1}{\widetilde K_{0}},
\end{aligned}
\end{equation*}
which concludes the proof.
\end{proof}

\medskip

According to Propositions~\ref{theo:P1} and~\ref{theo:P2}, the exact solution ${\mathbf Y}(t)$ and the \typo{RK--BUG} solution are bounded on the finite time-interval $0 \leq t \leq T$. \typo{We define}
\begin{equation}
\V_r := \Big\{\typo{{\mathbf Z}} \in \M_r \,\Big|\, \norm{\typo{{\mathbf Z}}-{\mathbf Y}_0}_F \leq \max\{R_1,R_2\} \Big\},
\end{equation}
where
\begin{equation*}
\begin{aligned}
R_1 &= \int_0^T e^{L(T-t)} \norm{{\mathbf F}(t,{\mathbf Y}_0)}_F \, \mathrm{d}t, \\
R_2 &= \max_{1 \leq i \leq s} \Bigg\{ ( 1+K_{i0} h_0 ) \frac{ e^{T \widetilde K_{0}}-1}{\widetilde K_{0}} \Bigg( \sum_{j=1}^{s} \widetilde K_{j} \Bigg) + h_0 \sum_{j=1}^{i-1} K_{ij} \Bigg\} \times \sup_{t \in \ff{0}{T}} \norm{{\mathbf F}(t,{\mathbf Y}_0)}_F.
\end{aligned}
\end{equation*}
The subset $\V_r \subseteq \M_r$ is a neighbourhood of the \modif{initial condition ${\mathbf Y}_0$} \typo{that} contains \typo{both the exact solution and its RK--BUG approximation} for all $t \in \ff{0}{T}$ and $h \leq h_0$. \typo{Consequently, as} ${\mathbf F}(t,{\mathbf X})$ is continuous in time and Lipschitz continuous in ${\mathbf X}$, ${\mathbf F}$ is bounded on the compact set $\ff{0}{T} \times \V_r$, and there exists a constant $B \geq 0$ such that
\begin{equation}
\label{eq:def_B}
\typo{B := \sup_{t \in \ff{0}{T}} \sup_{{\mathbf Z} \in \V_r} \norm{{\mathbf F}(t,{\mathbf Z})}_F.}
\end{equation}
\modif{Lastly, note that $\V_r$, and therefore $B$, are independent of~$h$}.

\medskip

\begin{lemma}
\label{theo:L4}
\typo{Suppose that Assumption~\ref{theo:A1} holds. Then,} the exact flow is $e^{Lt}$-Lipschitz continuous:
\begin{equation}
\norm{{\boldsymbol\Phi}^{t}_{\mathbf F}(\typo{{\mathbf Z}}_1)-{\boldsymbol\Phi}^{t}_{\mathbf F}(\typo{{\mathbf Z}}_2)}_F \leq e^{Lt} \norm{\typo{{\mathbf Z}}_1-\typo{{\mathbf Z}}_2}_F,
\end{equation}
for all $t \in \ff{0}{T}$ and $\typo{{\mathbf Z}}_1, \typo{{\mathbf Z}}_2 \in \R^{n \times m}$.
\end{lemma}
\begin{proof}
From Assumption \ref{theo:A1}, we obtain the differential inequality
\begin{equation*}
\begin{aligned}
\frac{\diff}{\diff t} \norm{{\boldsymbol\Phi}^{t}_{\mathbf F}(\typo{{\mathbf Z}}_1)-{\boldsymbol\Phi}^{t}_{\mathbf F}(\typo{{\mathbf Z}}_2)}_F^2 &= 2 \, \dott{{\boldsymbol\Phi}^{t}_{\mathbf F}(\typo{{\mathbf Z}}_1)-{\boldsymbol\Phi}^{t}_{\mathbf F}(\typo{{\mathbf Z}}_2)}{{\mathbf F}(t,{\boldsymbol\Phi}^{t}_{\mathbf F}(\typo{{\mathbf Z}}_1))-{\mathbf F}(t,{\boldsymbol\Phi}^{t}_{\mathbf F}(\typo{{\mathbf Z}}_2))}_F \\
&\leq 2 L \norm{{\boldsymbol\Phi}^{t}_{\mathbf F}(\typo{{\mathbf Z}}_1)-{\boldsymbol\Phi}^{t}_{\mathbf F}(\typo{{\mathbf Z}}_2)}_F^2.
\end{aligned}
\end{equation*}
\typo{Then}, according to Gr\"onwall's inequality, the exact flow \typo{satisfies}
\begin{equation*}
\norm{{\boldsymbol\Phi}^{t}_{\mathbf F}(\typo{{\mathbf Z}}_1)-{\boldsymbol\Phi}^{t}_{\mathbf F}(\typo{{\mathbf Z}}_2)}^2_F \leq e^{2Lt} \norm{\typo{{\mathbf Z}}_1-\typo{{\mathbf Z}}_2}^2_F,
\end{equation*}
which concludes the proof.
\end{proof}

\subsection{\modif{Intermediate error estimates}}

\modif{We now provide three error estimates that will be essential to establish the high-order convergence of the RK--BUG integrator. The first one, presented in Proposition~\ref{theo:P3}, concerns the projection error and shows that the Galerkin projection onto the augmented bases is more accurate than the tangent-space projection. This property will allow us to adapt the convergence analysis of PRK methods~\cite{kieri2019projection} to the present RK--BUG integrator. The second estimate, given in equation~\eqref{eq:def_gamma}, introduces a bound on the truncation error. Compared to~\cite{kieri2019projection}, this estimate is new and will allow us to derive error bounds with an improved order of convergence. Finally, the third estimate corresponds to the standard local error of high-order Runge--Kutta methods.}

\medskip

\begin{proposition}
\label{theo:P3}
\typo{The Galerkin projection onto the augmented bases provides a more accurate approximation of ${\mathbf F}_{k i}$ than the tangent-space projection. Specifically, \modif{for all $k \in \{0,\ldots,N-1\}$}, the projection error satisfies}
\begin{equation}
\label{eq:p3_1}
\norm{{\mathbf F}_{k i}-\widehat{\mathbf U}_{k+1}\widehat{\mathbf U}_{k+1}^\top{\mathbf F}_{k i}\widehat{\mathbf V}_{k+1}\widehat{\mathbf V}_{k+1}^\top}_F \leq \norm{{\mathbf F}_{k i}-\Po_{{\mathbf Y}_{k i}}\big({\mathbf F}_{k i}\big)}_F,
\end{equation}
for all $i \in \{1,\ldots,s\}$ such that $b_i \neq 0$, \typo{and}
\begin{equation}
\label{eq:p3_2}
\modif{\norm{{\mathbf F}_{k j}-\widehat{\mathbf U}_{k i}\widehat{\mathbf U}_{k i}^\top{\mathbf F}_{k j}\widehat{\mathbf V}_{k i}\widehat{\mathbf V}_{k i}^\top}_F \leq \norm{{\mathbf F}_{k j}-\Po_{{\mathbf Y}_{k j}}\big({\mathbf F}_{k j}\big)}_F,}
\end{equation}
for all \modif{$i \in \{2,\ldots,s\}$ and $j \in \{1,\ldots,i-1\}$ such that $a_{i j} \neq 0$}.
\end{proposition}
\begin{proof}
We prove equation~\eqref{eq:p3_1}. The proof of~\eqref{eq:p3_2} follows from the same arguments and is therefore omitted. \typo{Let the thin SVD of the tangent-space projection be
\(
\Po_{{\mathbf Y}_{k i}}\big({\mathbf F}_{k i}\big)
= {\mathbf \Phi}\,{\boldsymbol \Sigma}\,{\mathbf \Psi}^\top,
\)
with $\bar r \leq \hat r$ its rank. By construction of the augmented bases (for $b_i \neq 0$), the column and row spaces of $\Po_{{\mathbf Y}_{k i}}\big({\mathbf F}_{k i}\big)$ 
are in the spans of the augmented bases:
\begin{equation*}
\operatorname{range}({\mathbf \Phi})
  \subseteq \operatorname{span}(\widehat{\mathbf U}_{k+1}), 
\qquad
\operatorname{range}({\mathbf \Psi})
  \subseteq \operatorname{span}(\widehat{\mathbf V}_{k+1}).
\end{equation*}}
Consequently, the Galerkin projection onto the augmented bases \typo{satisfies}
\begin{equation*}
\begin{aligned}
\norm{{\mathbf F}_{k i}-\Po_{{\mathbf Y}_{k i}}\big({\mathbf F}_{k i}\big)}_F
&= \norm{{\mathbf F}_{k i}-{\mathbf \Phi}{\boldsymbol \Sigma}{\mathbf \Psi}^\top}_F \\
&\ge \min_{\overline{\boldsymbol \Sigma} \in \R^{\bar r \times \bar r}}
   \norm{{\mathbf F}_{k i}-{\mathbf \Phi}\overline{\boldsymbol \Sigma}{\mathbf \Psi}^\top}_F \\
&\ge \min_{\widehat{\boldsymbol \Sigma} \in \R^{\hat r \times \hat r}}
   \norm{{\mathbf F}_{k i}-\widehat{\mathbf U}_{k+1}\widehat{\boldsymbol \Sigma}\widehat{\mathbf V}_{k+1}^\top}_F \\
&= \norm{{\mathbf F}_{k i}
         -\widehat{\mathbf U}_{k+1}\widehat{\mathbf U}_{k+1}^\top
          {\mathbf F}_{k i}
          \widehat{\mathbf V}_{k+1}\widehat{\mathbf V}_{k+1}^\top}_F,
\end{aligned}
\end{equation*}
which concludes the proof.
\end{proof}

\medskip

Since ${\mathbf F}$ is bounded on $\ff{0}{T} \times \V$, \typo{its} orthogonal projection onto the tangent space is also bounded on $\ff{0}{T} \times \V_r$, and there exists $\varepsilon_r \geq 0$ such that the tangent-space projection error is bounded on $\ff{0}{T} \times \V_r$. \typo{We define} 
\begin{equation}
\label{eq:def_epsilon}
\varepsilon_r := \sup_{t \in \ff{0}{T}} \; \sup_{\typo{{\mathbf Z}} \in \V_r} \norm{{\mathbf F}(t,\typo{{\mathbf Z}}) - \Po_\typo{{\mathbf Z}}\big({\mathbf F}(t,\typo{{\mathbf Z}})\big)}_F.
\end{equation}
\typo{According to Proposition~\ref{theo:P3}}, the projection error of the \typo{RK--BUG} integrator 
\typo{satisfies}, \modif{for all $k \in \{0,\ldots,N-1\}$},
\begin{align*}
\norm{{\mathbf F}_{k i} -\widehat{\mathbf U}_{k+1}\widehat{\mathbf U}_{k+1}^\top {\mathbf F}_{k i}\widehat{\mathbf V}_{k+1}\widehat{\mathbf V}_{k+1}^\top}_F &\leq \varepsilon_r, \quad \modif{i \in \I_b}, \\
\norm{{\mathbf F}_{k j} -\widehat{\mathbf U}_{k i}\widehat{\mathbf U}_{k i}^\top {\mathbf F}_{k j} \widehat{\mathbf V}_{k i}\widehat{\mathbf V}_{k i}^\top}_F &\leq \varepsilon_r, \quad \modif{i \in \{2,\ldots,s\},\ j \in \I^{(i)}_a},
\end{align*}
\modif{where $\I_b = \{\, i \in \{1,\ldots,s\} \mid b_{i}\neq 0 \,\}$ and $\I^{(i)}_a = \{\, j \in \{1,\ldots,i-1\} \mid a_{i j}\neq 0 \,\}$}. \typo{Finally}, note that the tangent-space projection error, and therefore $\varepsilon_r$, vanish when $r = \min\{n,m\}$.

\medskip

\begin{proposition}
\label{theo:P4}
\review{Suppose that Assumption~\ref{theo:A1} holds. Then, for $h \leq h_0$ and a fixed rank~$r$, the truncation error satisfies, for all $k \in \{0,\ldots,N-1\}$,
\begin{align}
\norm{\widehat{\mathbf Y}_{k+1} - \Ro_{\M_r}\bigl( \widehat{\mathbf Y}_{k+1} \bigr)}_F &\leq h \overline\gamma_r, \label{eq:p4_1} \\
\norm{\widehat{\mathbf Y}_{k,i+1} - \Ro_{\M_r}\bigl( \widehat{\mathbf Y}_{k,i+1} \bigr)}_F &\leq h \overline\gamma_r, \quad i \in \{1,\ldots,s-1\}, \label{eq:p4_2}
\end{align}
where $\overline\gamma_r \geq 0$ is independent of $h$ and vanishes when $r = \min\{n,m\}$.}
\end{proposition}
\begin{proof}
We prove equation~\eqref{eq:p4_1}. The proof of~\eqref{eq:p4_2} follows from the same arguments and is therefore omitted. Let $\sigma_j(\typo{{\mathbf Z}})$ denote the $j$-th largest singular value of $\typo{{\mathbf Z}} \in \R^{n \times m}$. According to the Eckart--Young theorem~\cite{eckart1936approximation}, the squared truncation error is given by
\begin{equation}
\label{eq:l4_3}
\norm{\widehat{\mathbf Y}_{k+1} - \typo{\Ro_{\M_r}}\bigl( \widehat{\mathbf Y}_{k+1} \bigr)}_F^2 
  = \sum_{j=r+1}^{\min\{n,m,2rs\}} \sigma_j^2(\widehat{\mathbf Y}_{k+1}),
\end{equation}
since the rank of $\widehat{\mathbf Y}_{k+1}$ is at most $2rs$, and \typo{thus} $\sigma_j(\widehat{\mathbf Y}_{k+1}) = 0$ for $j > 2rs$. \typo{Moreover}, as $\widehat{\mathbf Y}_{k+1} 
  = {\mathbf Y}_{k} 
    + h \typo{\sum_{i=1}^s b_i \widetilde{\mathbf F}_{k i}}$ with $
  \widetilde{\mathbf F}_{k i} 
    = \widehat{\mathbf U}_{k+1}\widehat{\mathbf U}_{k+1}^\top
      {\mathbf F}_{k i}
      \widehat{\mathbf V}_{k+1}\widehat{\mathbf V}_{k+1}^\top$,
the singular values of $\widehat{\mathbf Y}_{k+1}$ are bounded, according to Theorem~3.3.16 in~\cite{horn1991topics}, by
\begin{equation*}
\sigma_{i+j-1}(\widehat{\mathbf Y}_{k+1}) 
  \leq \sigma_i({\mathbf Y}_{k}) 
       + h \, \sigma_j\!\left(\typo{\sum_{l=1}^s b_{l} \widetilde{\mathbf F}_{k l}}\right)
\end{equation*}
for all \typo{$i,j \in \{1,\ldots,\min\{n,m\}\}$} such that \typo{$i+j-1 \leq \min\{n,m\}$}. In particular, for $i = r+1$ and \typo{any} $j \ge 1$, it follows that
\begin{equation}
\label{eq:l4_4}
\sigma_{r+j}(\widehat{\mathbf Y}_{k+1}) 
  \leq h \, \sigma_j\!\left(\typo{\sum_{l=1}^s b_{l} \widetilde{\mathbf F}_{k l}}\right),
\end{equation}
since the rank of ${\mathbf Y}_{k}$ is at most $r$. \typo{Furthermore}, from equation~\eqref{eq:def_B}, \typo{we obtain}
\begin{equation*}
\begin{aligned}
\sum_{j=1}^{\min\{n,m,2rs\}} \sigma_j^2\!\left(\typo{\sum_{l=1}^s b_{l} \widetilde{\mathbf F}_{k l}}\right) 
 &= \Big\| \typo{\sum_{l=1}^s b_{l} \widetilde{\mathbf F}_{k l}} \Big\|_F^2 \\
 &\leq \Big( \typo{\sum_{l=1}^s |b_{l}| \, \|\widetilde{\mathbf F}_{k l} \|_F} \Big)^2 \\
 &\leq \Big( \typo{\sum_{l=1}^s |b_{l}| \, \|{\mathbf F}_{k l} \|_F} \Big)^2 \\
 &\leq \Big( \typo{\sum_{l=1}^s |b_{l}| \, B} \Big)^2 
  = C_b^2 B^2,
\end{aligned}
\end{equation*}
and consequently,
\begin{equation}
\label{eq:l4_5}
\sigma_{j}^2\!\left(\typo{\sum_{l=1}^s b_{l} \widetilde{\mathbf F}_{k l}}\right) \leq C_b^2 B^2, \qquad j \geq 1,
\end{equation}
where $C_b = \sum_{i=1}^s |b_i|$. Finally, combining~\eqref{eq:l4_3}--\eqref{eq:l4_5} yields
\begin{equation*}
\begin{aligned}
\norm{\widehat{\mathbf Y}_{k+1} - \typo{\Ro_{\M_r}}\bigl( \widehat{\mathbf Y}_{k+1} \bigr)}_F^2 &= \sum_{j=1}^{\min\{n,m,2rs\}-r} \sigma_{r+j}^2(\widehat{\mathbf Y}_{k+1}) \\
 &\leq h^2 \sum_{j=1}^{\min\{n,m,2rs\}-r} \sigma_j^2\!\left(\typo{\sum_{l=1}^s b_{l} \widetilde{\mathbf F}_{k l}}\right) \\
 &\leq h^2 (\min\{n,m,2rs\}-r) C_b^2 B^2.
\end{aligned}
\end{equation*}
\modif{The desired result follows directly, since $(\min\{n,m,2rs\}-r) C_b^2 B^2$ is independent of $h$ and vanishes when $r$ is full rank, as $(\min\{n,m,2rs\}-r)$ becomes zero while $C_b^2 B^2$ is finite.}
\end{proof}

\medskip

\review{According to Proposition~\ref{theo:P4}, there exists $\gamma_{r} \leq \overline\gamma_r$ such that~\eqref{eq:p4_1}--\eqref{eq:p4_2} hold. In the following, we define $\gamma_r$ as the smallest constant for which these bounds are valid:
\begin{equation}
\label{eq:def_gamma}
\begin{aligned}
\gamma_r := \frac{1}{h} 
\max\Big\{ &
\max_{k \in \{0,\ldots,N-1\}} \ \sup_{{\mathbf Z}_k \in \V_r} 
  \norm{\widehat{\mathbf Z}_{k+1} - \Ro_{\M_r}\bigl(\widehat{\mathbf Z}_{k+1}\bigr)}_F , \\
&\max_{\substack{k \in \{0,\ldots,N-1\} \\ i \in \{1,\ldots,s-1\}}} \ \sup_{{\mathbf Z}_{k} \in \V_r} 
  \norm{\widehat{\mathbf Z}_{k,i+1} - \Ro_{\M_r}\bigl(\widehat{\mathbf Z}_{k,i+1}\bigr)}_F 
\Big\}.
\end{aligned}
\end{equation}
where $\widehat{\mathbf Z}_{k+1}$ and $\widehat{\mathbf Z}_{k,i+1}$ denote the augmented RK--BUG solutions obtained by starting from ${\mathbf Z}_{k}$.}

\medskip

\begin{remark}
\modif{The error term $\gamma_r$ quantifies the truncation error and is generally unrelated to $\varepsilon_r$, which measures the projection error onto the tangent space.}
\end{remark}

\medskip

\begin{proposition}
\label{theo:P5}
\typo{Let $(a_{ij},b_j,c_i)$ define an explicit Runge--Kutta method~\eqref{eq:rk} of order~$p$}, and suppose that Assumption~\ref{theo:A2} holds. Then, for $h \leq h_0$, the local error of the Runge--Kutta method \typo{satisfies}
\begin{equation}
\norm{\typo{{\mathbf X}}_{k+1}-{\boldsymbol\Phi}^{h}_{\mathbf F}(\typo{{\mathbf X}}_k)}_F \leq C_L h^{p+1},
\end{equation}
where the constant $C_L>0$ is independent of $h$.
\end{proposition}
\begin{proof}
See Theorem~3.1 \modif{in} Chapter~2 \modif{of} \cite{harrier1993solving}.
\end{proof}

\subsection{\modif{Main convergence theorems}}

We now establish the local and global error bounds \typo{of the RK--BUG integrator} in Theorems~\ref{theo:T1} and~\ref{theo:T2}, respectively. These bounds show that the proposed integrator retains the order of convergence of the \typo{underlying} Runge--Kutta method until the error reaches a plateau corresponding to the low-rank truncation error, which vanishes as the rank becomes full. Compared to~\cite{kieri2019projection}, we obtain error bounds with an improved order \typo{of convergence}, thanks to the additional error term~$\gamma_r$. \typo{However}, similar bounds can be derived for \typo{PRK} methods by \typo{introducing~$\gamma_r$} and adapting the convergence analysis accordingly.

\medskip

\begin{lemma}
\label{theo:L3}
Suppose that Assumption~\ref{theo:A1} holds, \typo{and let} ${\mathbf Z}_{k i}$ denote the \typo{(unprojected) Runge--Kutta} solution obtained by starting from ${\mathbf Z}_{k} = {\mathbf Y}_{k} \in \M_r$. Then, for $h \leq h_0$ and \modif{a fixed rank~$r$}, \typo{the RK--BUG solution ${\mathbf Y}_{k i}$ satisfies}
\begin{equation}
\norm{{\mathbf Y}_{k i}-{\mathbf Z}_{k i}}_F \leq C_{i} h \left(\varepsilon_r + \gamma_r\right),
\end{equation}
\review{on each subinterval $\ff{t_k}{t_{k+1}} \subseteq \ff{0}{T}$ and for all stages $i \in \{1,\ldots,s\}$. Here, $C_{i}\geq 0$ is a constant independent of~$h$ and~$r$, $\varepsilon_r$ bounds the projection error (see~\eqref{eq:def_epsilon}), and $\gamma_r$ is defined so that $h\gamma_r$ bounds the truncation error (see~\eqref{eq:def_gamma}).}
\end{lemma}
\begin{proof}
We proceed by induction \typo{on the stage index~$i$}. For $i=1$, the statement is trivial with $C_1=0$, since ${\mathbf Z}_{k 1} = {\mathbf Z}_{k} = {\mathbf Y}_{k} = {\mathbf Y}_{k 1}$. \typo{Assume now that the result holds for all stages $j < i$, with $i \in \{2,\ldots,s\}$}. Then, the local error can be \typo{decomposed into two contributions}:
\begin{equation}
\label{eq:l3_1}
\norm{{\mathbf Y}_{k i}-{\mathbf Z}_{k i}}_F 
  \leq \norm{{\mathbf Y}_{k i}-\widehat{\mathbf Y}_{k i}}_F 
     + \norm{\widehat{\mathbf Y}_{k i}- {\mathbf Z}_{k i}}_F.
\end{equation}
\typo{According to Proposition~\ref{theo:P4}, the first contribution satisfies}
\begin{equation}
\label{eq:l3_2}
\norm{\widehat{\mathbf Y}_{k i}-{\mathbf Y}_{k i}}_F \leq h \gamma_r.
\end{equation}
\typo{For the second contribution, the induction hypothesis yields}
\begin{equation}
\label{eq:l3_3}
\begin{aligned}
\norm{\widehat{\mathbf Y}_{k i}-{\mathbf Z}_{k i}}_F 
&\leq h \sum_{j=1}^{i-1} |a_{i j}| \,
   \norm{\widehat{\mathbf U}_{k i}\widehat{\mathbf U}_{k i}^\top
   {\mathbf F}(t_{k j},{\mathbf Y}_{k j})
   \widehat{\mathbf V}_{k i}\widehat{\mathbf V}_{k i}^\top
   -{\mathbf F}(t_{k j},{\mathbf Z}_{k j})}_F \\
&= h \sum_{j \in \typo{\I^{(i)}_a}} |a_{i j}| \,
   \norm{\widehat{\mathbf U}_{k i}\widehat{\mathbf U}_{k i}^\top
   {\mathbf F}(t_{k j},{\mathbf Y}_{k j})
   \widehat{\mathbf V}_{k i}\widehat{\mathbf V}_{k i}^\top
   -{\mathbf F}(t_{k j},{\mathbf Z}_{k j})}_F \\
&\leq h \sum_{j \in \typo{\I^{(i)}_a}} |a_{i j}| \,
   \Bigl( \norm{\widehat{\mathbf U}_{k i}\widehat{\mathbf U}_{k i}^\top
   {\mathbf F}(t_{k j},{\mathbf Y}_{k j})
   \widehat{\mathbf V}_{k i}\widehat{\mathbf V}_{k i}^\top
   -{\mathbf F}(t_{k j},{\mathbf Y}_{k j})}_F \\
&\hphantom{\leq h \sum_{j \in \typo{\I^{(i)}_a}} |a_{i j}| \, \Bigl(}
   + \norm{{\mathbf F}(t_{k j},{\mathbf Y}_{k j})
   -{\mathbf F}(t_{k j},{\mathbf Z}_{k j})}_F \Bigr) \\
&\leq h \sum_{j \in \typo{\I^{(i)}_a}} |a_{i j}| \,
   \left( \varepsilon_r + L \norm{{\mathbf Y}_{k j}-{\mathbf Z}_{k j}}_F \right) \\
&\leq \Big( \sum_{j \in \typo{\I^{(i)}_a}} |a_{i j}| (1+L C_{j} h) \Big) h \varepsilon_r 
   + \Big( \sum_{j \in \typo{\I^{(i)}_a}} |a_{i j}| L C_{j} \Big) h^2 \gamma_r,
\end{aligned}
\end{equation}
where $\typo{\I^{(i)}_a} = \{\, j \in \{1,\ldots,i-1\} \mid a_{i j}\neq 0 \,\}$. Finally, combining~\eqref{eq:l3_1}--\eqref{eq:l3_3}, \typo{we obtain}
\begin{equation*}
\begin{aligned}
\norm{{\mathbf Y}_{k i}-{\mathbf Z}_{k i}}_F 
&\leq \Big( \sum_{j \in \typo{\I^{(i)}_a}} |a_{i j}| (1+L C_{j} h) \Big) h \varepsilon_r 
   + \Big( 1 + \sum_{j \in \typo{\I^{(i)}_a}} |a_{i j}| L C_{j} h \Big) h \gamma_r \\
&\leq C_i h (\varepsilon_r + \gamma_r),
\end{aligned}
\end{equation*}
where
\begin{equation*}
\begin{aligned}
C_1 &= 0, \\
C_i &= \max\Big\{
   \sum_{j=1}^{i-1} |a_{i j}| (1+L C_{j} h_0),
   \, 1+ \sum_{j=1}^{i-1} |a_{i j}| L C_{j} h_0
\Big\}, \quad i=2,\ldots,s.
\end{aligned}
\end{equation*}
\end{proof}

\medskip

\begin{theorem}[\modif{Local error bound}]
\label{theo:T1}
Suppose that Assumptions~\ref{theo:A1} and~\ref{theo:A2} hold. Then, for $h \leq h_0$ \modif{and a fixed rank~$r$}, the local error of the \typo{RK--BUG} integrator \typo{satisfies}
\begin{equation}
\norm{{\mathbf Y}_{k+1}-{\boldsymbol\Phi}^{h}_{\mathbf F}({\mathbf Y}_k)}_F \leq C h \left( \varepsilon_r + \gamma_r + h^{p} \right),
\end{equation}
\review{on each subinterval $\ff{t_k}{t_{k+1}} \subseteq \ff{0}{T}$. Here, $C>0$ is a constant independent of~$h$ and~$r$, $\varepsilon_r$ bounds the projection error (see~\eqref{eq:def_epsilon}), and $\gamma_r$ is defined so that $h\gamma_r$ bounds the truncation error (see~\eqref{eq:def_gamma}).}
\end{theorem}
\begin{proof}
Let ${\mathbf Z}_{k+1}$ denote the \typo{Runge--Kutta solution obtained by starting from} ${\mathbf Z}_{k} = {\mathbf Y}_{k} \in \M_r$. The local error can be \typo{decomposed into three contributions:}
\begin{equation}
\label{eq:t1_1}
\norm{{\mathbf Y}_{k+1}-{\boldsymbol\Phi}^{h}_{\mathbf F}({\mathbf Y}_k)}_F  
  \leq \norm{{\mathbf Y}_{k+1}- \widehat{\mathbf Y}_{k+1}}_F 
     + \norm{\widehat{\mathbf Y}_{k+1}- {\mathbf Z}_{k+1}}_F 
     + \norm{{\mathbf Z}_{k+1}-{\boldsymbol\Phi}^{h}_{\mathbf F}({\mathbf Y}_k)}_F.
\end{equation}
\typo{According to Proposition~\ref{theo:P4}, the first contribution satisfies}
\begin{equation}
\label{eq:t1_2}
\norm{\widehat{\mathbf Y}_{k+1}-{\mathbf Y}_{k+1}}_F \leq h \gamma_r.
\end{equation}
\typo{For the second contribution, using Lemma~\ref{theo:L3}, we obtain}
\begin{equation}
\label{eq:t1_3}
\begin{aligned}
\norm{\widehat{\mathbf Y}_{k+1}-{\mathbf Z}_{k+1}}_F 
&\leq h \sum_{i=1}^{s} |b_{i}| \,
   \norm{\widehat{\mathbf U}_{k+1}\widehat{\mathbf U}_{k+1}^\top
   {\mathbf F}(t_{k i},{\mathbf Y}_{k i})
   \widehat{\mathbf V}_{k+1}\widehat{\mathbf V}_{k+1}^\top
   -{\mathbf F}(t_{k i},{\mathbf Z}_{k i})}_F \\
&= h \sum_{i \in \I_b} |b_{i}| \,
   \norm{\widehat{\mathbf U}_{k+1}\widehat{\mathbf U}_{k+1}^\top
   {\mathbf F}(t_{k i},{\mathbf Y}_{k i})
   \widehat{\mathbf V}_{k+1}\widehat{\mathbf V}_{k+1}^\top
   -{\mathbf F}(t_{k i},{\mathbf Z}_{k i})}_F \\
&\leq h \sum_{i \in \I_b} |b_{i}| \,
   \Bigl( \norm{\widehat{\mathbf U}_{k+1}\widehat{\mathbf U}_{k+1}^\top
   {\mathbf F}(t_{k i},{\mathbf Y}_{k i})
   \widehat{\mathbf V}_{k+1}\widehat{\mathbf V}_{k+1}^\top
   -{\mathbf F}(t_{k i},{\mathbf Y}_{k i})}_F \\
&\hphantom{\leq h \sum_{i \in \I_b} |b_{i}| \, \Bigl(}
   + \norm{{\mathbf F}(t_{k i},{\mathbf Y}_{k i})
   -{\mathbf F}(t_{k i},{\mathbf Z}_{k i})}_F \Bigr) \\
&\leq h \sum_{i \in \I_b} |b_{i}| 
   \left(\varepsilon_r + L \norm{{\mathbf Y}_{k i}-{\mathbf Z}_{k i}}_F\right) \\
&\leq \Big( \sum_{i \in \I_b} |b_{i}| (1+L C_{i} h) \Big) h \varepsilon_r 
   + \Big( \sum_{i \in \I_b} |b_{i}| L C_{i} \Big) h^2 \gamma_r,
\end{aligned}
\end{equation}
where $\I_b = \{\, i \in \{1,\ldots,s\} \mid b_{i}\neq 0 \,\}$. \typo{The last contribution is bounded, according to Proposition~\ref{theo:P5}, by}
\begin{equation}
\label{eq:t1_4}
\norm{{\mathbf Z}_{k+1}-{\boldsymbol\Phi}^{h}_{\mathbf F}({\mathbf Y}_k)}_F 
  = \norm{{\mathbf Z}_{k+1}-{\boldsymbol\Phi}^{h}_{\mathbf F}({\mathbf Z}_k)}_F 
  \leq C_L h^{p+1}.
\end{equation}
Finally, combining~\eqref{eq:t1_1}--\eqref{eq:t1_4} yields the desired result:
\begin{equation*}
\begin{aligned}
\norm{{\mathbf Y}_{k+1}-{\boldsymbol\Phi}^{h}_{\mathbf F}({\mathbf Y}_k)}_F 
&\leq \Big( \sum_{i \in \I_b} |b_{i}| (1+L C_{i} h) \Big) h \varepsilon_r 
   + \Big( 1 + \sum_{i \in \I_b} |b_{i}| L C_{i} h \Big) h \gamma_r 
   + C_L h^{p+1} \\
&\leq C h (\varepsilon_r + \gamma_r + h^p),
\end{aligned}
\end{equation*}
where
\begin{equation*}
C = \max\Big\{ 
   \sum_{i=1}^s |b_{i}| (1+L C_{i} h_0), \,
   1 + \sum_{i=1}^s |b_{i}| L C_{i} h_0, \,
   C_L 
\Big\}.
\end{equation*}
\end{proof}

\medskip

\begin{theorem}[\modif{Global error bound}]
\label{theo:T2}
Suppose that Assumptions~\ref{theo:A1} and~\ref{theo:A2} hold. Then, for $h \leq h_0$ \modif{and a fixed rank~$r$}, the global error of the \typo{RK--BUG} integrator \typo{satisfies}
\begin{equation}
\label{eq:global_error}
\norm{\typo{{\mathbf X}}(t_N)-{\mathbf Y}_N}_F \leq C' \left( \delta_r + \varepsilon_r + \gamma_r + h^{p} \right),
\end{equation}
\review{on the finite time-interval $0 \leq t_N = N h \leq T$. Here, ${\mathbf X}(t_N)$ denotes the exact solution of~\eqref{eq:ode}, and ${\mathbf Y}_N$ its RK--BUG approximation. The constant $C'>0$ depends on $T$, $L$, $C_L$ (see Proposition~\ref{theo:P5}), $h_0$, $C_a = \sum_{i,j=1}^s |a_{ij}|$, and $C_b = \sum_{i=1}^s |b_i|$, but not on~$h$ or~$r$. Finally, $\delta_r := \norm{{\mathbf X}_0-{\mathbf Y}_0}_F$ is the initial error, $\varepsilon_r$ bounds the projection error (see~\eqref{eq:def_epsilon}), and $\gamma_r$ is defined so that $h\gamma_r$ bounds the truncation error (see~\eqref{eq:def_gamma}); all these terms vanish when $r = \min\{n,m\}$.}
\end{theorem}
\begin{proof}
\typo{The result follows} from the local error \typo{bound} of Theorem~\ref{theo:T1} and the standard argument of Lady Windermere's fan, \typo{which describes the accumulation of local errors} along the exact flow. \typo{Specifically}, the global error can be \typo{written} as the telescoping sum
\begin{equation}
\label{eq:t2_1}
{\mathbf Y}_{N}-{\boldsymbol\Phi}^{Nh}_{\mathbf F}(\typo{{\mathbf X}}_0) 
  = \sum_{k=1}^N 
    \Big( 
      {\boldsymbol\Phi}^{(N-k)h}_{\mathbf F}({\mathbf Y}_{k}) 
      - {\boldsymbol\Phi}^{(N-k+1)h}_{\mathbf F}({\mathbf Y}_{k-1}) 
    \Big)
    + {\boldsymbol\Phi}^{Nh}_{\mathbf F}({\mathbf Y}_{0}) 
      - {\boldsymbol\Phi}^{Nh}_{\mathbf F}(\typo{{\mathbf X}}_{0}).
\end{equation}
\typo{The different contributions are} bounded, according to Lemma \ref{theo:L4} and Theorem \ref{theo:T1}, \typo{as follows:}
\begin{equation}
\label{eq:t2_2}
\norm{{\boldsymbol\Phi}^{Nh}_{\mathbf F}({\mathbf Y}_{0}) 
      - {\boldsymbol\Phi}^{Nh}_{\mathbf F}(\typo{{\mathbf X}}_{0})}_F
  \leq e^{LNh} \norm{{\mathbf Y}_{0} - \typo{{\mathbf X}}_{0}}_F 
  = e^{LNh} \delta_r,
\end{equation}
\begin{equation}
\label{eq:t2_3}
\begin{aligned}
\norm{{\boldsymbol\Phi}^{(N-k)h}_{\mathbf F}({\mathbf Y}_{k}) 
      - {\boldsymbol\Phi}^{(N-k)h}_{\mathbf F}({\boldsymbol\Phi}^{h}_{\mathbf F}({\mathbf Y}_{k-1}))}_F 
&\leq e^{L(N-k)h} 
      \norm{{\mathbf Y}_{k}-{\boldsymbol\Phi}^{h}_{\mathbf F}({\mathbf Y}_{k-1})}_F \\
&\leq e^{L(N-k)h} \, C h (\varepsilon_r + \gamma_r + h^p).
\end{aligned}
\end{equation}
\typo{Combining~\eqref{eq:t2_1}--\eqref{eq:t2_3} then yields}
\begin{equation*}
\norm{{\mathbf Y}_{N}-{\boldsymbol\Phi}^{Nh}_{\mathbf F}(\typo{{\mathbf X}}_0)}_F 
  \leq e^{LNh} \delta_r 
  + C (\varepsilon_r + \gamma_r + h^p) 
    \sum_{k=1}^N e^{L(N-k)h} h.
\end{equation*}
Finally, the Riemann sum \typo{satisfies}
\begin{equation*}
\sum_{k=1}^N h e^{L(N-k)h} 
  \leq \int_0^{Nh} e^{L(Nh-t)} \diff t 
  = \frac{e^{LNh}-1}{L},
\end{equation*}
and using $Nh \leq T$, \typo{we obtain}
\begin{equation*}
\begin{aligned}
\norm{{\mathbf Y}_{N}-{\boldsymbol\Phi}^{Nh}_{\mathbf F}(\typo{{\mathbf X}}_0)}_F 
&\leq e^{LT} \delta_r 
  + C \frac{e^{LT}-1}{L} (\varepsilon_r + \gamma_r + h^p) \\
&\leq \max\!\left\{ e^{LT},\, C \frac{e^{LT}-1}{L} \right\} 
   (\delta_r + \varepsilon_r + \gamma_r + h^p),
\end{aligned}
\end{equation*}
which concludes the proof.
\end{proof}

\medskip

\begin{remark}
\typo{When} $\max\{ \delta_r, \varepsilon_r, \gamma_r \} \ll h^p$, the \typo{RK--BUG} integrator \typo{retains} the order of convergence of the \typo{underlying} Runge--Kutta method. Otherwise, the global error is dominated by the initial, projection, and truncation errors, \typo{which} vanish as the rank becomes full.
\end{remark}

\subsection{\review{Rank-adaptive strategy}}
\label{sec:rank_adaptive}

\review{In the RK--BUG integrator, an adaptive rank can be used instead of a fixed one. A common approach is to truncate the augmented solution so that the truncation error remains below a prescribed tolerance. According to Theorem~\ref{theo:T2}, this tolerance must depend on the step size~$h$ to maintain the convergence order~$p$. Specifically, the rank~$r$ should be chosen such that the error terms $\delta_r$, $\varepsilon_r$, and $\gamma_r$ are proportional to~$h^p$. Hence, Theorem~\ref{theo:T2} provides a practical way to design rank-adaptive strategies. One such strategy, which can be combined with more sophisticated approaches, is to select the rank as the smallest integer such that the augmented solution satisfies, at each intermediate stage $i \in \{1,\ldots,s-1\}$,
\begin{equation*}
\norm{\widehat{\mathbf Y}_{k,i+1} - \Ro_{\M_r}\bigl( \widehat{\mathbf Y}_{k,i+1} \bigr)}_F \leq \alpha h^{p+1},
\end{equation*}
and, at the final stage,
\begin{equation*}
\norm{\widehat{\mathbf Y}_{k+1} - \Ro_{\M_r}\bigl( \widehat{\mathbf Y}_{k+1} \bigr)}_F \leq \alpha h^{p+1},
\end{equation*}
where $\alpha \geq 0$ is a tunable constant. In practice, we also employ a relative tolerance~$\beta$ (taken as~$\beta = 10^{-14}$ in the numerical experiments), since $h^{p+1}$ can become very small for high orders~$p$. Moreover, we choose an initial rank~$r_0$ that is not too small and enforce the subsequent ranks to always remain above~$r_0$ in order to prevent an early truncation of modes associated with singular values that are initially small but may become important later in time. In summary, for each time-step $k \in \{0,\ldots,N-1\}$, the augmented solution is truncated so that $r \ge r_0$ and
\begin{equation*}
\begin{aligned}
\norm{\widehat{\mathbf Y}_{k,i+1} - \Ro_{\M_r}\bigl( \widehat{\mathbf Y}_{k,i+1} \bigr)}_F &\leq \max\!\left\{ \alpha h^{p+1}, \; \beta \norm{\widehat{\mathbf Y}_{k,i+1}}_F \right\},
\quad i \in \{1,\ldots,s-1\}, \\
\norm{\widehat{\mathbf Y}_{k+1} - \Ro_{\M_r}\bigl( \widehat{\mathbf Y}_{k+1} \bigr)}_F &\leq \max\!\left\{ \alpha h^{p+1}, \; \beta \norm{\widehat{\mathbf Y}_{k+1}}_F \right\}.
\end{aligned}
\end{equation*}}

\section{Numerical experiments}
\label{sec:4}

\review{In this section, we assess the performance of the proposed RK--BUG integrator through several numerical experiments. The objectives of these experiments are fourfold:
\begin{enumerate}
\item validate the high-order convergence of the RK--BUG integrator for different explicit Runge--Kutta schemes, including second-order (midpoint, Heun), third-order, and fourth-order methods (see Appendix~\ref{sec:A1} for the corresponding Butcher tableaux);
\item compare the accuracy of the RK--BUG integrator with existing dynamical low-rank integrators, such as the midpoint BUG integrator~\cite{ceruti2024robust} (using its first variant, which is more accurate than the second one) and PRK methods~\cite{kieri2019projection};
\item illustrate the rank-adaptive strategy;
\item verify that the conservative RK--BUG variant preserves physical invariants, such as total mass and momentum.
\end{enumerate}
The first three aspects are investigated on three benchmark problems (the Allen--Cahn, Lyapunov, and discrete nonlinear Schr\"odinger equations) taken from~\cite{kieri2019projection,lam2024randomized}, while the conservative variant is evaluated on the Vlasov--Poisson equations. For the adaptive-rank experiments, the parameter $\alpha$ is chosen empirically to obtain a small adaptive rank while preserving high-order convergence. The initial rank $r_0$ is determined from the fixed-rank experiments to achieve a comparable error level. Throughout this section}, the accuracy is measured with respect to a reference solution $\typo{\mathbf{X}}_k$ computed using the Runge--Kutta--Fehlberg method, and the error is defined as
\begin{equation*}
\mathrm{Error} = \max_{0 \leq k \leq N} \| \typo{\mathbf{X}}_k - \mathbf{Y}_k \|_F.
\end{equation*}

\subsection{Allen-Cahn equation}

We first consider the Allen--Cahn equation:
\begin{equation*}
\typo{\dot{\mathbf X}} =  \theta ({\mathbf L} \typo{{\mathbf X}} + \typo{{\mathbf X}} {\mathbf L}) + \typo{{\mathbf X}} - \typo{{\mathbf X}} \odot \typo{{\mathbf X}} \odot \typo{{\mathbf X}}, \qquad \typo{{\mathbf X}}(0) = \typo{{\mathbf X}}_0,
\end{equation*}
where $\typo{{\mathbf X}}(t) \in \R^{n \times n}$, ${\mathbf L} = \frac{n^2}{4 \pi^2}\typo{\mathrm{tridiag}}(1,-2,1) \in \R^{n \times n}$, $t \in \ff{0}{10}$, $\theta = 10^{-2}$, and $\odot$ stands for the Hadamard product. The domain $\ff{0}{2\pi}^2$ is discretized using \typo{$n \times n$ equidistant grid points $(x_i, y_j)$, with $n = 128$}, and the initial condition is given by 
\begin{equation*}
(\typo{{\mathbf X}}_0)_{ij} = \frac{\left[e^{-\tan^2(x_i)}+e^{-\tan^2(y_j)}\right]\sin(x_i)\sin(y_j)}{1+e^{|\csc(-x_i/2)|}+e^{|\csc(-y_j/2)|}}.
\end{equation*}

\begin{figure}
\begin{subfigure}[c]{0.5\linewidth}
\center
\includegraphics[width=6.2cm]{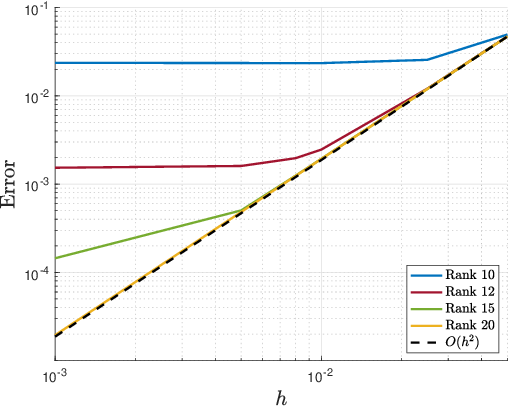}
\caption{RK--BUG (Midpoint)}
\end{subfigure}\hfill
\begin{subfigure}[c]{0.5\linewidth}
\center
\includegraphics[width=6.2cm]{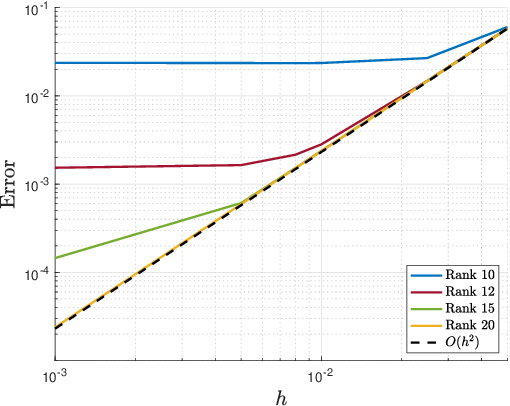}
\caption{RK--BUG (Heun)}
\end{subfigure}\hfill
\begin{subfigure}[c]{0.5\linewidth}
\center
\includegraphics[width=6.2cm]{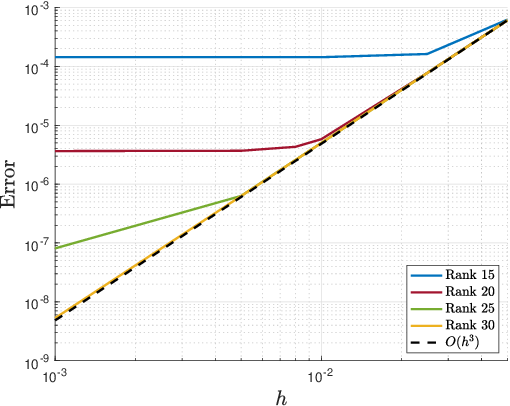}
\caption{RK--BUG (SSP33)}
\end{subfigure}\hfill
\begin{subfigure}[c]{0.5\linewidth}
\center
\includegraphics[width=6.2cm]{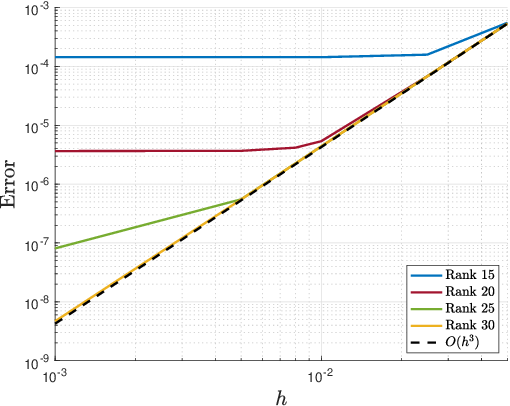}
\caption{RK--BUG (Heun3)}
\end{subfigure}\hfill
\begin{subfigure}[c]{1\linewidth}
\center
\includegraphics[width=6.2cm]{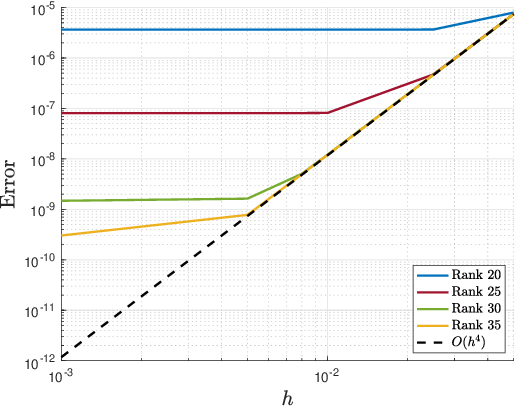}
\caption{RK--BUG (RK4)}
\end{subfigure}\hfill
\caption{\label{fig:allen-cahn_conv}Convergence error of high-order RK--BUG integrators for the Allen--Cahn equation. Dashed lines show reference slopes~$h^2, h^3, h^4$.}
\end{figure}

\noindent\textbf{\typo{High-order convergence.}}
\typo{Figure~\ref{fig:allen-cahn_conv}} presents the \typo{convergence} error of \typo{RK--BUG integrators with respect to} the step size $h$ \typo{for different} ranks $r$. \typo{The RK--BUG integrators} achieve second-, third-, and fourth-order convergence, \typo{depending on the underlying Runge--Kutta scheme}, until the error reaches a plateau corresponding to the low-rank truncation error. This plateau decreases as the rank increases, \typo{down} to a limit \typo{around} $10^{-9}$ \typo{due to} the accumulation of roundoff error.
\medskip

\begin{figure}
\begin{subfigure}[c]{0.5\linewidth}
\center
\includegraphics[width=6.2cm]{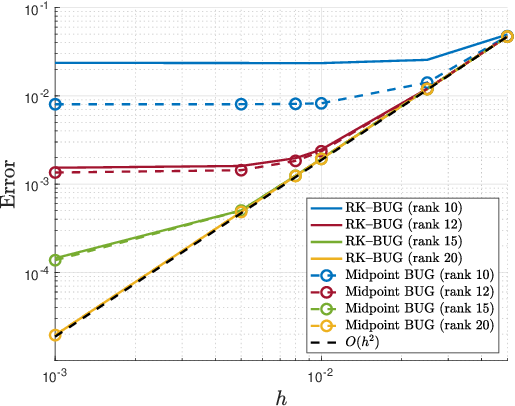}
\caption{RK--BUG (Midpoint) vs Midpoint BUG}
\end{subfigure}\hfill
\begin{subfigure}[c]{0.5\linewidth}
\center
\includegraphics[width=6.2cm]{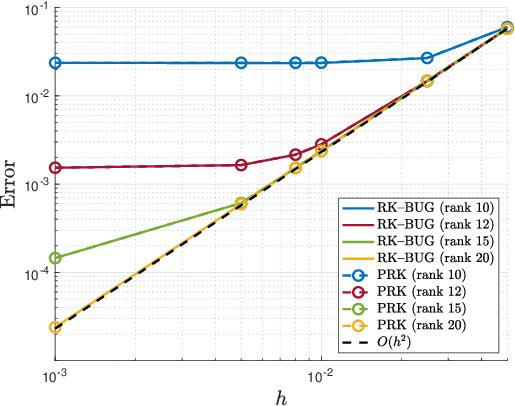}
\caption{RK--BUG (Heun) vs PRK (Heun)}
\end{subfigure}\hfill
\begin{subfigure}[c]{1\linewidth}
\center
\includegraphics[width=6.2cm]{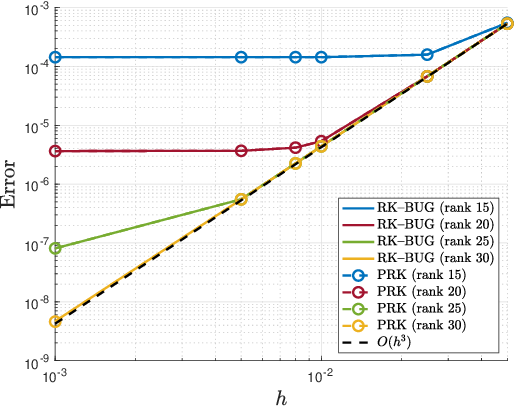}
\caption{RK--BUG (Heun3) vs PRK (Heun3)}
\end{subfigure}\hfill
\caption{\label{fig:allen-cahn_comp}Comparison of the RK--BUG integrator with other dynamical low-rank integrators for the Allen--Cahn equation.}
\end{figure}

\noindent\textbf{\typo{Comparison with existing methods.}}
\typo{Figure~\ref{fig:allen-cahn_comp} compares the RK--BUG integrator with the midpoint BUG and PRK methods.} \modif{For the midpoint scheme, the midpoint BUG integrator is slightly more accurate than the RK--BUG integrator for small ranks, while both approaches yield nearly identical errors as the rank increases.} \review{This difference can be explained by the fact that, unlike our RK--BUG (Midpoint) integrator, the midpoint BUG integrator does not truncate the intermediate solution, resulting in an augmented rank of at most $4r$ instead of $2r$. As a consequence, it is computationally more expensive but potentially more accurate.} \typo{For the other schemes, the RK--BUG and PRK methods exhibit the same accuracy.}
\medskip

\begin{figure}
\begin{subfigure}[c]{0.5\linewidth}
\center
\includegraphics[width=6.2cm]{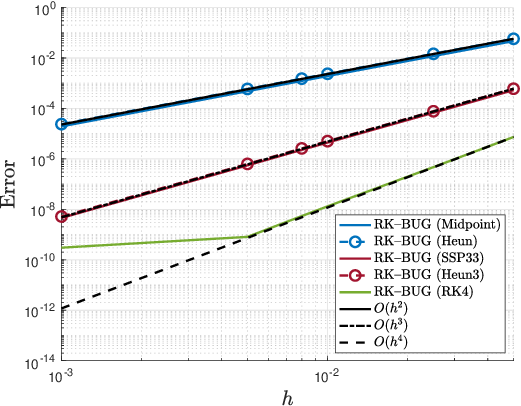}
\caption{Convergence error}
\end{subfigure}\hfill
\begin{subfigure}[c]{0.5\linewidth}
\center
\includegraphics[width=6.2cm]{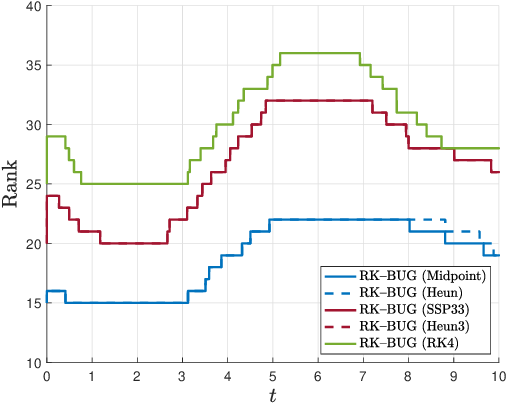}
\caption{Rank evolution for $h = 10^{-3}$}
\end{subfigure}\hfill
\caption{\label{fig:allen-cahn_adapt}Results of the rank-adaptive RK--BUG integrator for the Allen--Cahn equation.}
\end{figure}

\noindent\textbf{\review{Rank adaptivity.}}
\review{Figure~\ref{fig:allen-cahn_adapt} illustrates the performance of the rank-adaptive RK--BUG integrator. The parameter $\alpha$ is set to $10$, and the evolution of the rank $r$ is reported for $h=10^{-3}$, corresponding to the most demanding case. The results show that the adaptive RK--BUG integrator preserves the high-order convergence of the underlying Runge--Kutta scheme while using a smaller average rank compared to the fixed-rank RK--BUG integrator.}

\subsection{Lyapunov equation}

Then, we consider the \typo{(continuous-time)} Lyapunov equation:
\begin{equation*}
\typo{\dot{\mathbf X}} = {\mathbf L} \typo{{\mathbf X}} + \typo{{\mathbf X}} {\mathbf L} + \theta \frac{\mathbf C}{\norm{\mathbf C}_F}, \qquad \typo{{\mathbf X}}(0) = \typo{\mathbf{X}}_0,
\end{equation*}
where $\typo{{\mathbf X}}(t) \in \mathbb{R}^{n \times n}$, ${\mathbf L} = \frac{n^2}{4 \pi^2}\typo{\mathrm{tridiag}}(1,-2,1) \in \mathbb{R}^{n\times n}$, ${\mathbf C} \in \mathbb{R}^{n \times n}$, $t \in \ff{0}{10}$, and $\theta = 1$. The domain ${\ff{-\pi}{\pi}}^2$ is discretized using \typo{$n \times n$ equidistant grid points $(x_i, y_j)$, with $n = 128$}, and the initial condition and forcing term are defined as
\begin{equation*}
(\typo{{\mathbf X}}_0)_{ij} = \sin(x_i)\sin(y_j), \qquad ({\mathbf C})_{ij} = \sum_{l=1}^{11} 10^{-(l-1)} e^{-l(x_i^2+y_j^2)}.
\end{equation*}

\begin{figure}
\begin{subfigure}[c]{0.5\linewidth}
\center
\includegraphics[width=6.2cm]{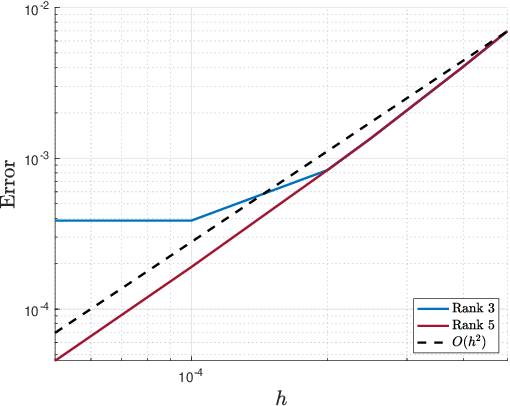}
\caption{RK--BUG (Midpoint)}
\end{subfigure}\hfill
\begin{subfigure}[c]{0.5\linewidth}
\center
\includegraphics[width=6.2cm]{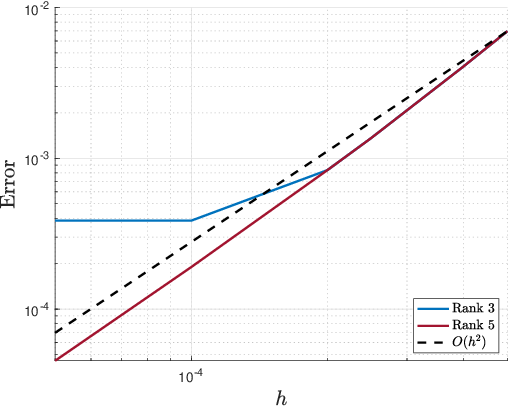}
\caption{RK--BUG (Heun)}
\end{subfigure}\hfill
\begin{subfigure}[c]{0.5\linewidth}
\center
\includegraphics[width=6.2cm]{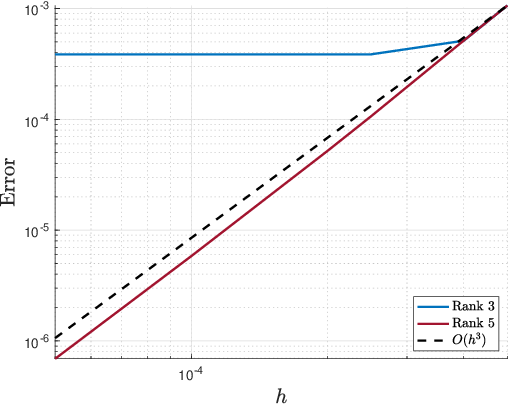}
\caption{RK--BUG (SSP33)}
\end{subfigure}\hfill
\begin{subfigure}[c]{0.5\linewidth}
\center
\includegraphics[width=6.2cm]{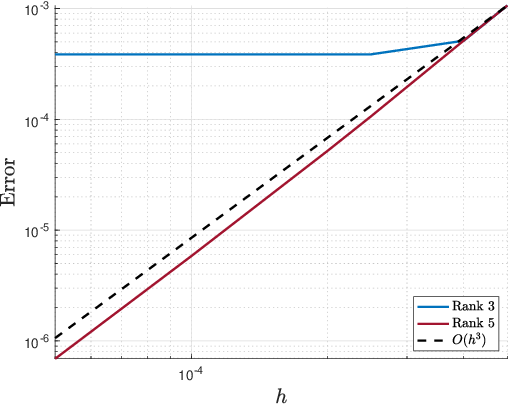}
\caption{RK--BUG (Heun3)}
\end{subfigure}\hfill
\begin{subfigure}[c]{1\linewidth}
\center
\includegraphics[width=6.2cm]{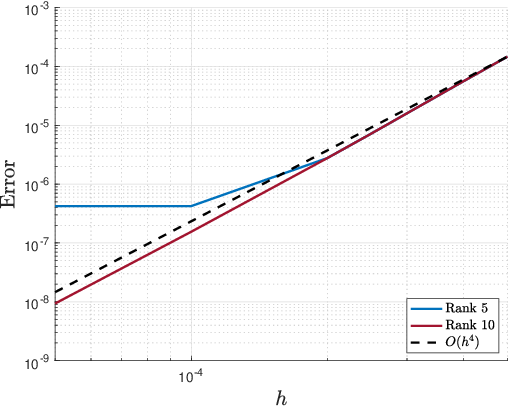}
\caption{RK--BUG (RK4)}
\end{subfigure}\hfill
\caption{\label{fig:lyapunov_conv}Convergence error of high-order RK--BUG integrators for the Lyapunov equation. Dashed lines show reference slopes~$h^2, h^3, h^4$.}
\end{figure}

\noindent\textbf{\typo{High-order convergence.}}
\typo{Figure~\ref{fig:lyapunov_conv} presents the convergence error of the RK--BUG integrator with respect to the step size $h$ for different ranks $r$.} \modif{The results confirm second-, third-, and fourth-order convergence, down to a plateau corresponding to the low-rank truncation error, which decreases as $r$ increases.}
\medskip

\begin{figure}
\begin{subfigure}[c]{0.5\linewidth}
\center
\includegraphics[width=6.2cm]{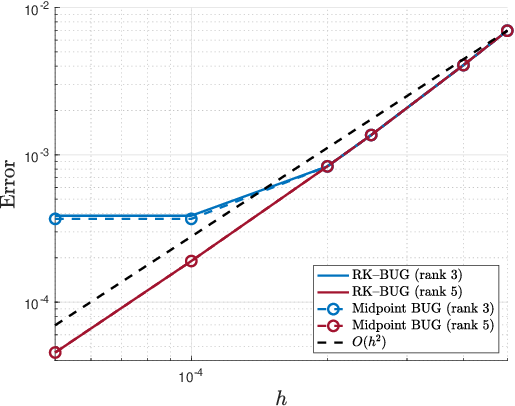}
\caption{RK--BUG (Midpoint) vs Midpoint BUG}
\end{subfigure}\hfill
\begin{subfigure}[c]{0.5\linewidth}
\center
\includegraphics[width=6.2cm]{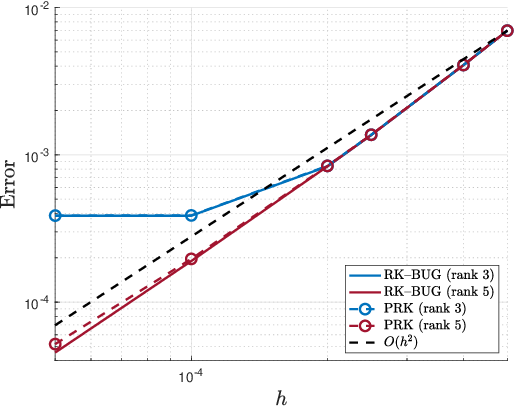}
\caption{RK--BUG (Heun) vs PRK (Heun)}
\end{subfigure}\hfill
\begin{subfigure}[c]{1\linewidth}
\center
\includegraphics[width=6.2cm]{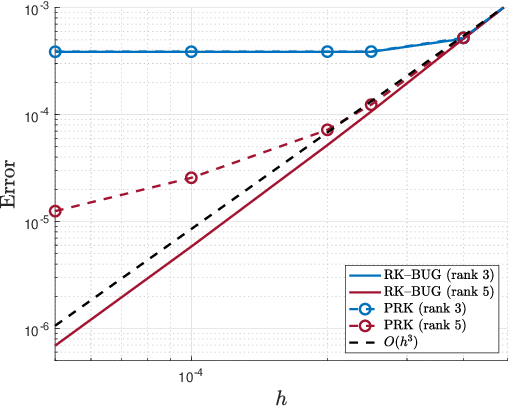}
\caption{RK--BUG (Heun3) vs PRK (Heun3)}
\end{subfigure}\hfill
\caption{\label{fig:lyapunov_comp}Comparison of the RK--BUG integrator with other dynamical low-rank integrators for the Lyapunov equation.}
\end{figure}

\noindent\textbf{\typo{Comparison with existing methods.}}
\typo{Figure~\ref{fig:lyapunov_comp} compares the RK--BUG integrator with the midpoint BUG and PRK methods.} \modif{For the second-order schemes, the errors of the midpoint BUG and PRK methods are almost identical to those of the corresponding RK--BUG integrators.} \typo{However, for the third-order scheme, the RK--BUG integrator is significantly more accurate than the PRK method. This difference can be explained by the fact that the Galerkin projection used in RK--BUG integrators provides a more accurate approximation of the discrete solution than the tangent-space projection employed in PRK methods (see Proposition~\ref{theo:P3}).}
\medskip

\begin{figure}
\begin{subfigure}[c]{0.5\linewidth}
\center
\includegraphics[width=6.2cm]{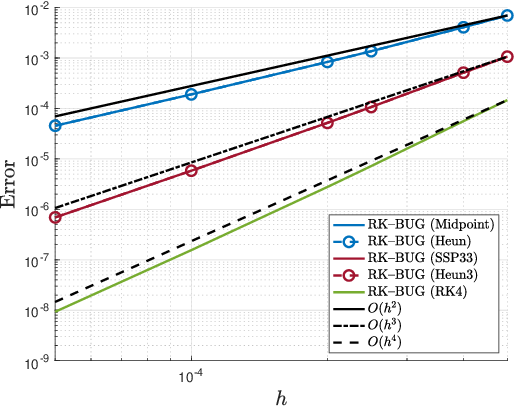}
\caption{Convergence error}
\end{subfigure}\hfill
\begin{subfigure}[c]{0.5\linewidth}
\center
\includegraphics[width=6.2cm]{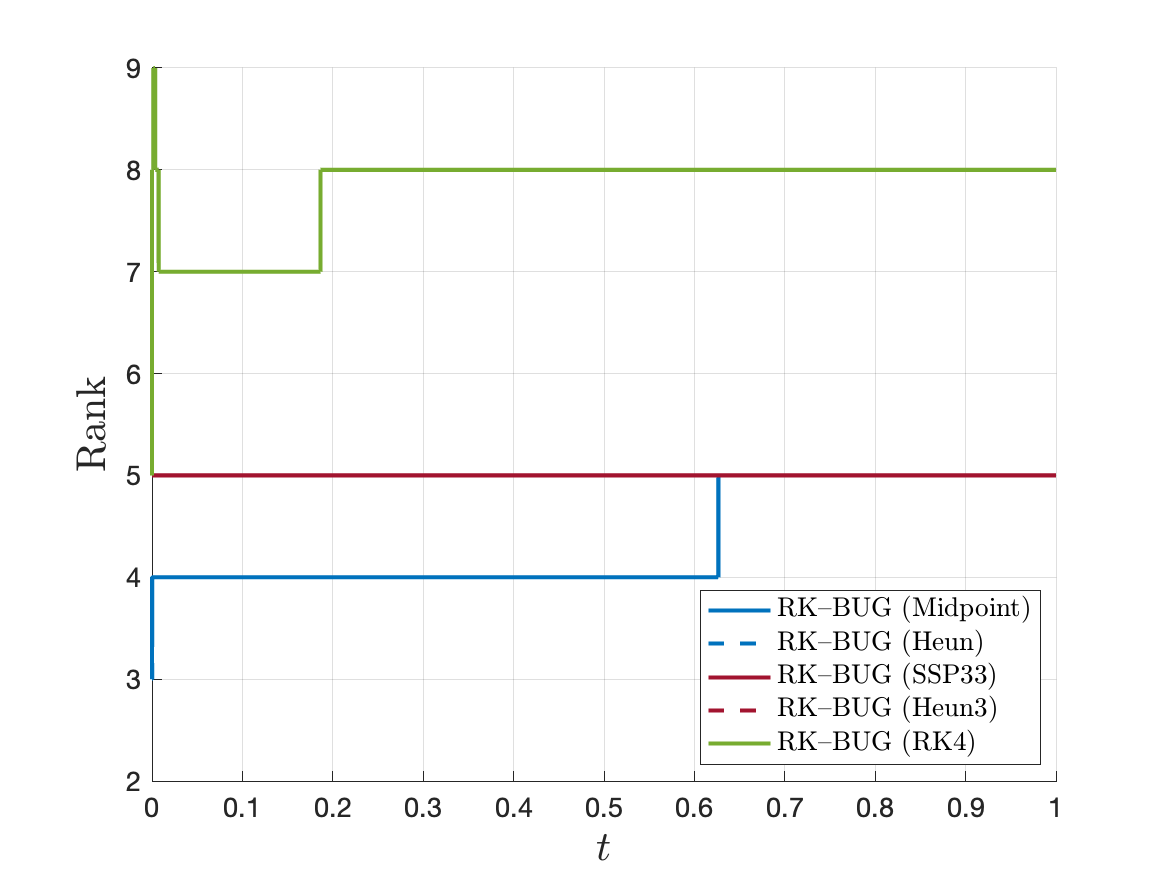}
\caption{Rank evolution}
\end{subfigure}\hfill
\caption{\label{fig:lyapunov_adapt}Results of the rank-adaptive RK--BUG integrator for the Lyapunov equation.}
\end{figure}

\noindent\textbf{\review{Rank adaptivity.}}
\review{Figure~\ref{fig:lyapunov_adapt} illustrates the performance of the rank-adaptive RK--BUG integrator. The parameter $\alpha$ is set to $10^5$ for $p=2$ and to $10^9$ for $p \geq 3$, and the evolution of the rank $r$ is reported for $h=5 \times 10^{-5}$, corresponding to the most demanding case. The results show that the adaptive RK--BUG integrator preserves the expected order of accuracy while using a smaller average rank than the fixed-rank RK--BUG integrator.}

\subsection{Discrete nonlinear Schr\"odinger equation}

We now consider the discrete nonlinear Schr\"odinger \typo{(DNLS)} equation:
\begin{equation*}
i \typo{\dot{\mathbf X}} = -\frac{1}{2} \left( {\mathbf D} \typo{{\mathbf X}} + \typo{{\mathbf X}} {\mathbf D} \right) - \theta \, |\typo{{\mathbf X}}|^2 \odot \typo{{\mathbf X}}, \qquad \typo{{\mathbf X}}(0) = \typo{{\mathbf X}}_0,
\end{equation*}
where $\typo{{\mathbf X}}(t) \in \mathbb{C}^{n \times n}$, ${\mathbf D} = \typo{\mathrm{tridiag}}(1,0,1) \in \mathbb{R}^{n \times n}$, \typo{$|{\mathbf X}|^2 := {\mathbf X} \odot \overline{\mathbf X}$ denotes the elementwise squared magnitude}, $t \in \ff{0}{5}$, $\theta = 0.3$, $n = 128$, and the initial condition is given by
\begin{equation*}
(\typo{{\mathbf X}}_0)_{jl} = \exp\!\left(-\frac{(j-60)^2}{100} - \frac{(l-50)^2}{100}\right) + \exp\!\left(-\frac{(j-50)^2}{100} - \frac{(l-40)^2}{100}\right).
\end{equation*}

\begin{figure}
\begin{subfigure}[c]{0.5\linewidth}
\center
\includegraphics[width=6.2cm]{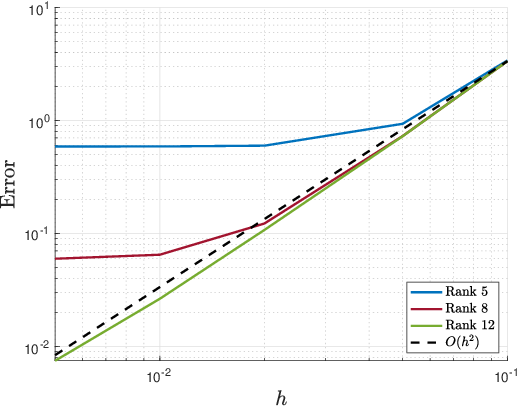}
\caption{RK--BUG (Midpoint)}
\end{subfigure}\hfill
\begin{subfigure}[c]{0.5\linewidth}
\center
\includegraphics[width=6.2cm]{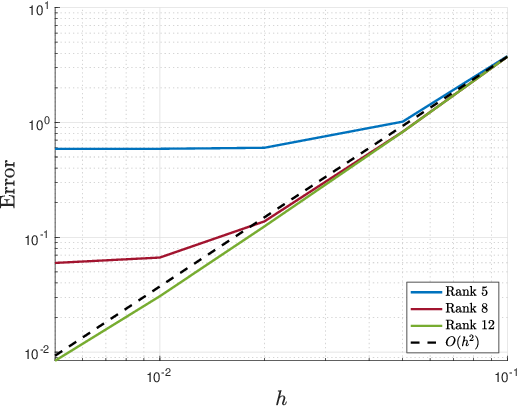}
\caption{RK--BUG (Heun)}
\end{subfigure}\hfill
\begin{subfigure}[c]{0.5\linewidth}
\center
\includegraphics[width=6.2cm]{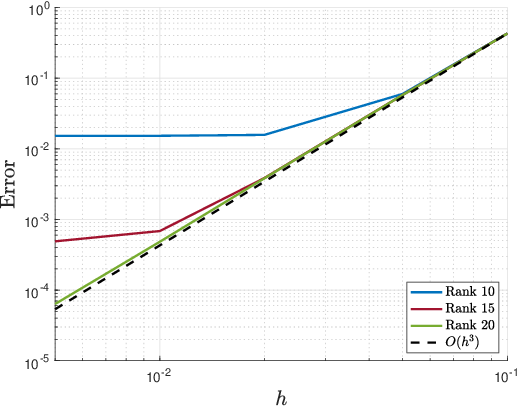}
\caption{RK--BUG (SSP33)}
\end{subfigure}\hfill
\begin{subfigure}[c]{0.5\linewidth}
\center
\includegraphics[width=6.2cm]{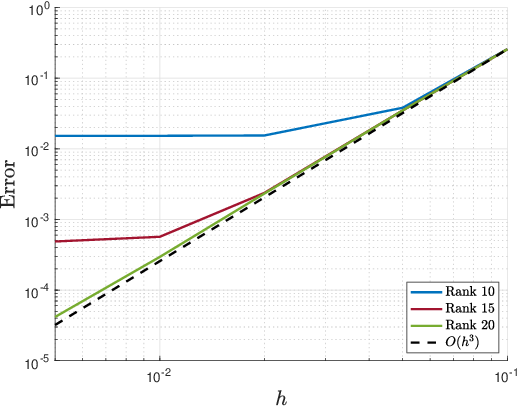}
\caption{RK--BUG (Heun3)}
\end{subfigure}\hfill
\begin{subfigure}[c]{1\linewidth}
\center
\includegraphics[width=6.2cm]{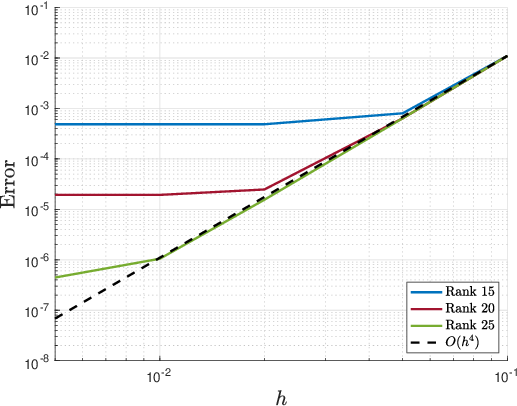}
\caption{RK--BUG (RK4)}
\end{subfigure}\hfill
\caption{\label{fig:dnls_conv}Convergence error of high-order RK--BUG integrators for the DNLS equation. Dashed lines show reference slopes~$h^2, h^3, h^4$.}
\end{figure}

\noindent\textbf{\typo{High-order convergence.}}
\typo{Figure~\ref{fig:dnls_conv} presents} the \typo{convergence} error of the \typo{RK--BUG} integrator \typo{with respect to} the step size $h$ \typo{for different} ranks $r$. \typo{The results confirm that the RK--BUG integrators achieve the expected second-, third-, and fourth-order convergence before reaching a plateau corresponding to the low-rank truncation error.}
\medskip

\begin{figure}
\begin{subfigure}[c]{0.5\linewidth}
\center
\includegraphics[width=6.2cm]{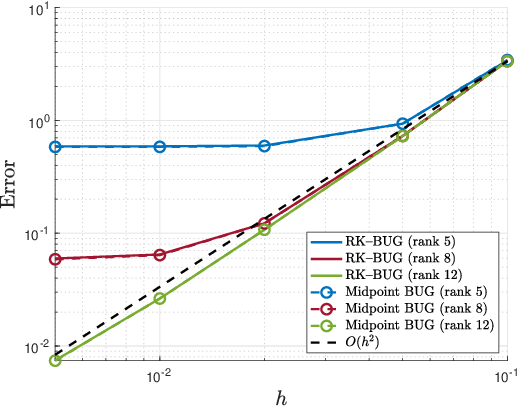}
\caption{RK--BUG (Midpoint) vs Midpoint BUG}
\end{subfigure}\hfill
\begin{subfigure}[c]{0.5\linewidth}
\center
\includegraphics[width=6.2cm]{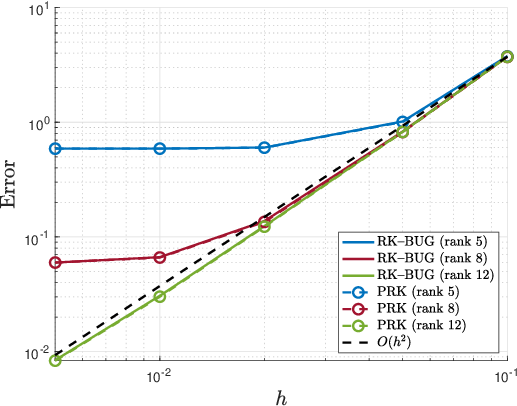}
\caption{RK--BUG (Heun) vs PRK (Heun)}
\end{subfigure}\hfill
\begin{subfigure}[c]{1\linewidth}
\center
\includegraphics[width=6.2cm]{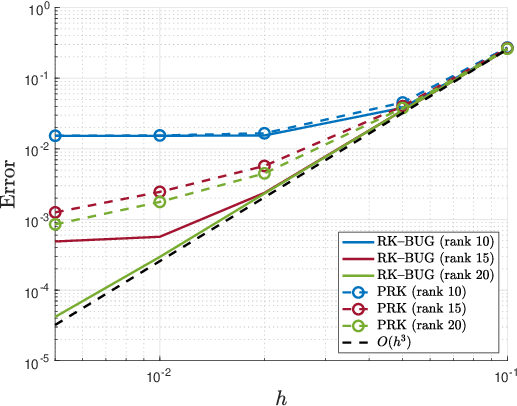}
\caption{RK--BUG (Heun3) vs PRK (Heun3)}
\end{subfigure}\hfill
\caption{\label{fig:dnls_comp}Comparison of the RK--BUG integrator with other dynamical low-rank integrators for the DNLS equation.}
\end{figure}

\noindent\textbf{\typo{Comparison with existing methods.}}
\typo{Figure~\ref{fig:dnls_comp} compares the RK--BUG integrator with the midpoint BUG and PRK methods. For the second-order schemes,} the errors of the midpoint BUG and \typo{PRK} methods are almost identical to those of the corresponding \typo{RK--BUG integrators}. However, for the \typo{third-order scheme}, the \typo{RK--BUG} integrator is significantly more accurate than the \typo{PRK} method.
\medskip

\begin{figure}
\begin{subfigure}[c]{0.5\linewidth}
\center
\includegraphics[width=6.2cm]{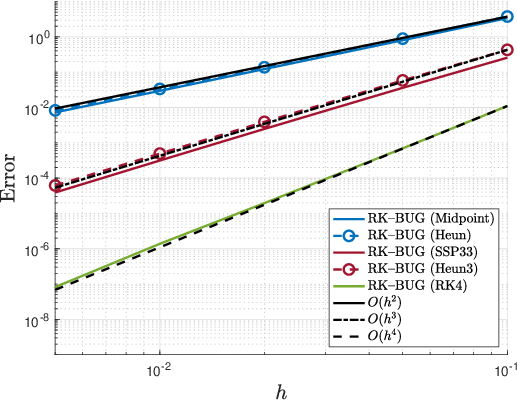}
\caption{Convergence error}
\end{subfigure}\hfill
\begin{subfigure}[c]{0.5\linewidth}
\center
\includegraphics[width=6.2cm]{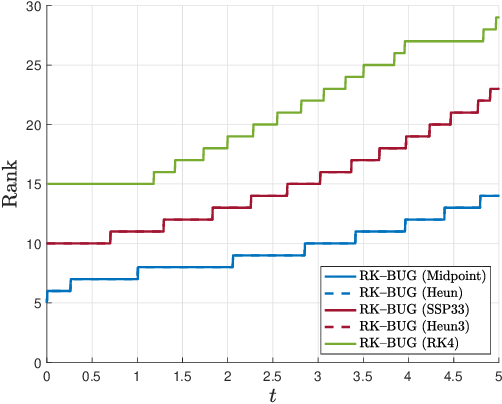}
\caption{Rank evolution}
\end{subfigure}\hfill
\caption{\label{fig:dnls_adapt}Results of the rank-adaptive RK--BUG integrator for the DNLS equation.}
\end{figure}

\noindent\textbf{\review{Rank adaptivity.}}
\review{Figure~\ref{fig:dnls_adapt} illustrates the performances of the rank-adaptive RK--BUG integrator. The parameter $\alpha$ is set to $10^2$, and the evolution of the rank $r$ is reported for $h=5 \times 10^{-3}$, corresponding to the most demanding case. The results show that the adaptive RK--BUG integrator preserves the high-order convergence of the underlying Runge--Kutta scheme while using a smaller average rank compared to the fixed-rank RK--BUG integrator.}

\subsection{\review{Vlasov--Poisson equations}}

\review{
Finally, we consider the one-dimensional in space and velocity (1D1V) Vlasov--Poisson equations with a constant background ion density:
\begin{equation*}
\begin{aligned}
\partial_t f(x,v,t) + v\,\partial_x f(x,v,t) - E(x,t)\,\partial_v f(x,v,t) = 0, \\
\partial_x E(x,t) = 1- \int_\R f(x,v,t) \diff v,
\end{aligned}
\end{equation*}
where $x \in \Omega_x \subset \R$ and $v \in \R$. The electric field $E(x,t)$ is computed here from the electron density $f(x,v,t)$ via the electrostatic potential $\phi(x,t)$, obtained by solving the Poisson equation:
\begin{equation*}
-\partial_x^2 \phi(x,t) = 1- \int_\R f(x,v,t) \diff v, \qquad
E(x,t) = -\,\partial_x \phi(x,t).
\end{equation*}
The Vlasov--Poisson equations admit three invariants: the total mass $N(t)$, momentum $J(t)$, and energy $\mathcal{E}(t)$,
\begin{equation*}
N(t) = \int_{\Omega_x} \rho(x,t)\,\mathrm{d}x, \qquad
J(t) = \int_{\Omega_x} j(x,t)\,\mathrm{d}x, \qquad
\mathcal{E}(t) = \int_{\Omega_x} e(x,t)\,\mathrm{d}x,
\end{equation*}
which are associated with the local quantities
\begin{equation*}
\begin{aligned}
\rho(x,t) &= \int_{\mathbb{R}} f(x,v,t)\,\mathrm{d}v,\\
j(x,t) &= \int_{\mathbb{R}} v\,f(x,v,t)\,\mathrm{d}v,\\
e(x,t) &= \frac{1}{2}\int_{\mathbb{R}} v^2 f(x,v,t)\,\mathrm{d}v + \frac{1}{2}E(x,t)^2.
\end{aligned}
\end{equation*}
We evaluate the conservative RK--BUG variant for the two-stream instability problem, defined by the initial condition
\begin{equation*}
f(x,v,0) = \left(1 + 10^{-3} \cos(0.2x)\right)
\frac{1}{2\sqrt{2\pi}}\left(\exp\!\left(-\frac{(v-2.4)^2}2\right) + \exp\!\left(-\frac{(v+2.4)^2}2\right)\right),
\end{equation*}
on the spatial domain $\Omega_x = [0,10\pi]$ and the truncated velocity domain $\Omega_v = [-9,9]$. Periodic boundary conditions are imposed in space, while artificial boundary conditions are applied in velocity due to the truncation of the infinite domain. In addition, we employ the smooth window function
\begin{equation*}
w(v)=
\begin{cases}
\exp\Bigl(\log(10^{-16})\,\bigl(\tfrac{v+7}{2}\bigr)^2\Bigr), & v<-7,\\[2pt]
1, & -7 \leq v \leq 7,\\[2pt]
\exp\Bigl(\log(10^{-16})\,\bigl(\tfrac{v-7}{2}\bigr)^2\Bigr), & v>7,
\end{cases}
\end{equation*}
when evaluating the charge density $\rho(x,t) = \int_{\Omega_v} f(x,v,t)\,w(v)\,\mathrm{d}v$ for the computation of the electric field, to ensure that the distribution function remains well-defined on the finite velocity domain and to suppress spurious reflections near the velocity boundaries.
}

\review{
The numerical schemes are chosen to preserve the conservation laws at the discrete level. The electric field is computed by solving the Poisson equation using the fast Fourier transform (FFT). Spatial and velocity derivatives are approximated by second-order summation-by-parts (SBP) upwind schemes \cite{mattsson2017diagonal}, which in particular satisfy the SBP property (see Appendix \ref{sec:A2} for details), and integrals are evaluated using the quadrature rule associated with the SBP norm. The computational domain $\Omega_x \times \Omega_v$ is discretized by a $128 \times 128$ uniform grid. Time integration is performed using various explicit Runge--Kutta methods with the fixed step size $h = 10^{-2}$.
}

\begin{figure}
\begin{subfigure}[c]{0.5\linewidth}
\center
\includegraphics[width=6.2cm]{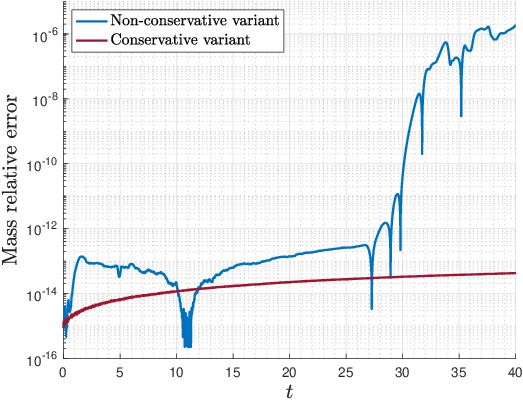}
\caption{RK--BUG (Midpoint)}
\end{subfigure}\hfill
\begin{subfigure}[c]{0.5\linewidth}
\center
\includegraphics[width=6.2cm]{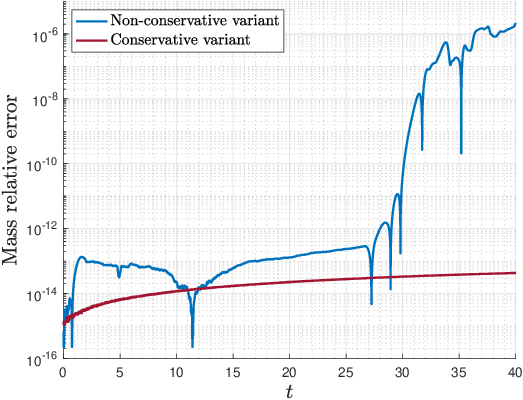}
\caption{RK--BUG (Heun)}
\end{subfigure}\hfill
\begin{subfigure}[c]{0.5\linewidth}
\center
\includegraphics[width=6.2cm]{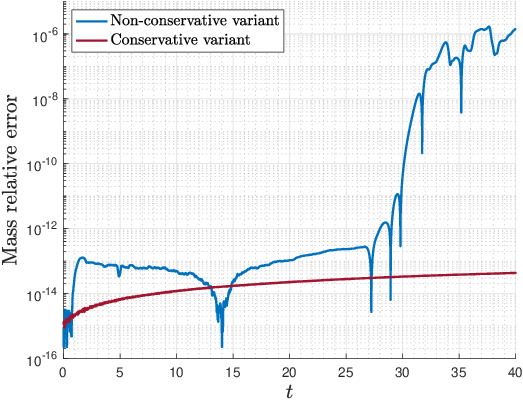}
\caption{RK--BUG (SSP33)}
\end{subfigure}\hfill
\begin{subfigure}[c]{0.5\linewidth}
\center
\includegraphics[width=6.2cm]{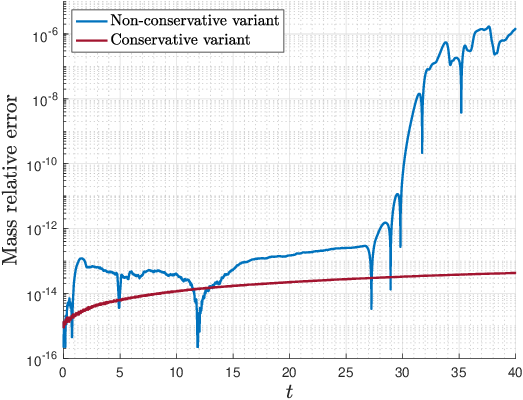}
\caption{RK--BUG (Heun3)}
\end{subfigure}\hfill
\begin{subfigure}[c]{1\linewidth}
\center
\includegraphics[width=6.2cm]{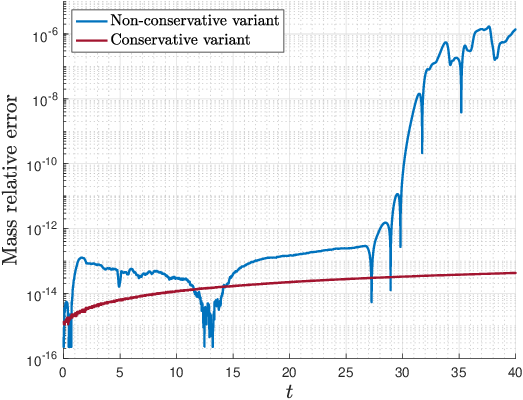}
\caption{RK--BUG (RK4)}
\end{subfigure}\hfill
\caption{\label{fig:vlasov_mass}Mass conservation errors of the conservative and non-conservative RK--BUG integrators for the Vlasov--Poisson equations.}
\end{figure}

\begin{figure}
\begin{subfigure}[c]{0.5\linewidth}
\center
\includegraphics[width=6.2cm]{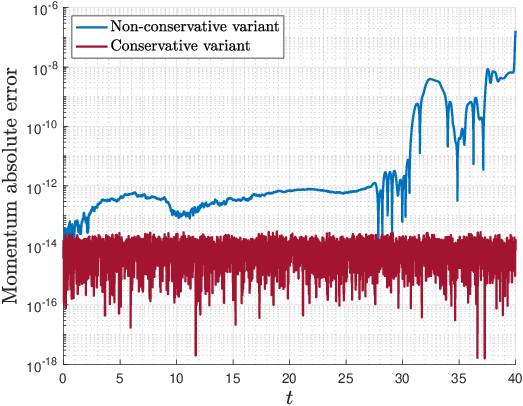}
\caption{RK--BUG (Midpoint)}
\end{subfigure}\hfill
\begin{subfigure}[c]{0.5\linewidth}
\center
\includegraphics[width=6.2cm]{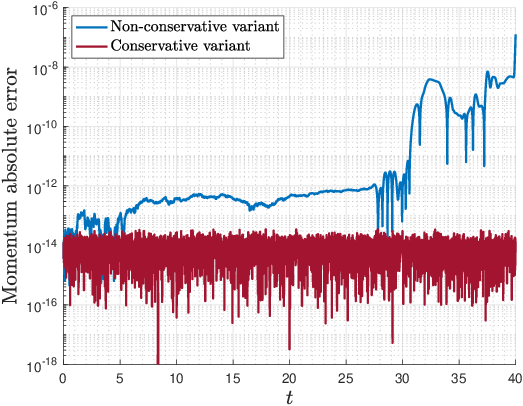}
\caption{RK--BUG (Heun)}
\end{subfigure}\hfill
\begin{subfigure}[c]{0.5\linewidth}
\center
\includegraphics[width=6.2cm]{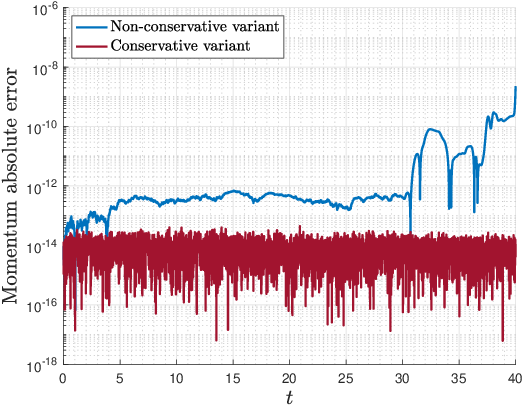}
\caption{RK--BUG (SSP33)}
\end{subfigure}\hfill
\begin{subfigure}[c]{0.5\linewidth}
\center
\includegraphics[width=6.2cm]{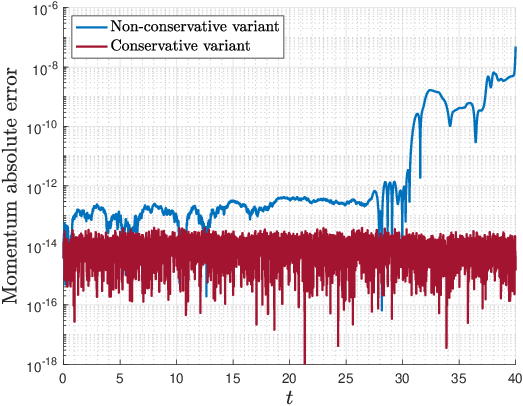}
\caption{RK--BUG (Heun3)}
\end{subfigure}\hfill
\begin{subfigure}[c]{1\linewidth}
\center
\includegraphics[width=6.2cm]{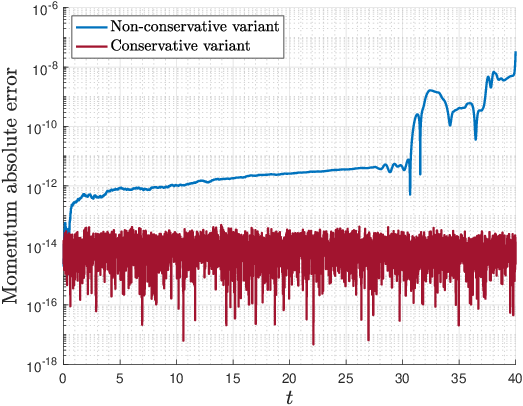}
\caption{RK--BUG (RK4)}
\end{subfigure}\hfill
\caption{\label{fig:vlasov_momentum}Momentum conservation errors of the conservative and non-conservative RK--BUG integrators for the Vlasov--Poisson equations.}
\end{figure}

\review{
We compare the standard RK--BUG integrator (using rank $r = 25$) with its conservative variant, constructed by enriching the basis with $\mathbf{W} = \mathbf{H}_v \, [\mathbf{1},\,\boldsymbol{v}] \in \mathbb{R}^{128 \times 2}$ (see Appendix~\ref{sec:A2} for the definition of $\mathbf{H}_v$), yielding $r_{\mathrm{cons}} = 2$ conservative modes in addition to the standard rank $r = 25$. The total energy is not included in this comparison, as explicit Runge--Kutta methods do not preserve it. Figures \ref{fig:vlasov_mass} and \ref{fig:vlasov_momentum} report the relative mass error and the absolute momentum error (the latter being zero). The overall behavior is very similar across all RK--BUG integrators. For the non-conservative variant, the error remains very low at early times but increases significantly as the simulation evolves. In contrast, the conservative RK--BUG integrator keeps the error negligible throughout the entire simulation.
}

\section{Conclusion}
\label{sec:5}

In this \typo{paper}, we have \typo{introduced} high-order BUG integrators based on explicit Runge--Kutta methods. \typo{These RK--BUG} integrators are robust \modif{with respect to small singular values, fully forward in time, and high-order accurate, while enabling conservation and rank adaptivity.} \typo{We have proved that RK--BUG integrators retain} the order of convergence of the \typo{underlying} Runge--Kutta method until the error reaches a plateau corresponding to the low-rank truncation error, which vanishes as the rank becomes full. \typo{The numerical experiments confirm the expected convergence orders $p=2$--$4$. Moreover, compared to existing dynamical low-rank integrators,} \modif{the RK--BUG integrator matches the accuracy of the midpoint BUG integrator for large ranks while requiring smaller augmented ranks}, \typo{and it significantly outperforms PRK methods for high orders ($p \geq 3$) and small step sizes.} \typo{In future work, we plan to extend the RK--BUG framework} to other classes of Runge--Kutta \typo{schemes}, such as exponential \typo{or implicit} Runge--Kutta methods.

\backmatter

\bmhead{Acknowledgements}

This work has been supported by the Swiss National Science Foundation under the Project n°200518 "Dynamical low rank methods for uncertainty quantification and data assimilation".

\section*{Declarations}

\bmhead{Data availability}

The datasets generated during and/or analyzed during the current study are available from the corresponding author on reasonable request.

\bmhead{Conflict of interest}

The authors declare that they have no conflict of interest.

\begin{appendices}

\section{List of Runge--Kutta methods}\label{sec:A1}

The Butcher tableaux associated with the Runge--Kutta methods used in this work are listed below.

\begin{itemize}
\item Forward Euler method
\medskip

\begin{tabular}{c | c c}
0 & 0 \\
1 & 1 \\
\hline
 & 1
\end{tabular}

\medskip
\item Explicit midpoint method \typo{(Midpoint)}
\medskip

\begin{tabular}{c | c c}
0 & 0 & 0 \\
1/2 & 1/2 & 0 \\
\hline
 & 0 & 1
\end{tabular}

\medskip
\item Heun's method \typo{(Heun)}
\medskip

\begin{tabular}{c | c c}
0 & 0 & 0 \\
1 & 1 & 0 \\
\hline
 & 1/2 & 1/2
\end{tabular}

\medskip
\item Third-order \typo{Strong Stability Preserving} Runge--Kutta method \typo{(SSP33)}
\medskip

\begin{tabular}{c | c c c}
0 & 0 & 0 & 0 \\
1 & 1 & 0 & 0 \\
1/2 & 1/4 & 1/4 & 0 \\
\hline
 & 1/6 & 1/6 & 2/3
\end{tabular}

\medskip
\item Heun's third-order method \typo{(Heun3)}
\medskip

\begin{tabular}{c | c c c}
0 & 0 & 0 & 0 \\
1/3 & 1/3 & 0 & 0 \\
2/3 & 0 & 2/3 & 0 \\
\hline
 & 1/4 & 0 & 3/4
\end{tabular}

\medskip
\item Classic fourth-order Runge--Kutta method \typo{(RK4)}
\medskip

\begin{tabular}{c | c c c c}
0 & 0 & 0 & 0 & 0 \\
1/2 & 1/2 & 0 & 0 & 0 \\
1/2 & 0 & 1/2 & 0 & 0 \\
1 & 0 & 0 & 1 & 0 \\
\hline
 & 1/6 & 1/3 & 1/3 & 1/6
\end{tabular}
\end{itemize}

\section{\review{Second-order SBP upwind schemes}}\label{sec:A2}

\review{
The derivatives in $x$ and $v$ are discretized using second-order SBP upwind finite-difference schemes:
\begin{equation*}
\begin{aligned}
v\,\partial_x f &\;\approx\; \max(v,0)\,\mathbf{D}_x^- f + \min(v,0)\,\mathbf{D}_x^+ f, \\[2pt]
- E\,\partial_v f &\;\approx\; \max(-E,0)\,\mathbf{D}_v^- f + \min(-E,0)\,\mathbf{D}_v^+ f.
\end{aligned}
\end{equation*}
These schemes are designed to mimic the integration-by-parts formula at the discrete level. In the $x$-direction, the operators $\mathbf{D}_x^+$ and $\mathbf{D}_x^-$ are defined as
\begin{equation*}
\mathbf{D}_x^+ = \frac{1}{2\Delta x}
\begin{bmatrix}
-3 & 4 & -1 & 0 & \cdots & 0 \\
0 & -3 & 4 & -1 & \ddots & \vdots \\
\vdots & \ddots & \ddots & \ddots & \ddots & 0 \\
0 & & \ddots & -3 & 4 & -1 \\
-1 & \ddots &  & \ddots & -3 & 4 \\
4 & -1 & 0 & \cdots & 0 & -3
\end{bmatrix},
\qquad
\mathbf{D}_x^- = -(\mathbf{D}_x^+)^{\!\top}.
\end{equation*}
They satisfy the discrete SBP property
\begin{equation*}
\mathbf{H}_x \mathbf{D}_x^+ + (\mathbf{H}_x \mathbf{D}_x^-)^{\!\top} = \mathbf{B}_x,
\end{equation*}
where $\mathbf{H}_x = \Delta x\,\mathbf{I}$ is the symmetric positive-definite matrix defining the SBP norm and, due to periodic boundary conditions, $\mathbf{B}_x = \mathbf{0}$. In the $v$-direction, the pair $(\mathbf{D}_v^+,\mathbf{D}_v^-)$ is obtained from
\begin{equation*}
\mathbf{Q}_v^+ = \frac{1}{4}
\begin{bmatrix}
-1 & 5 & -2 & 0 & \cdots & \cdots & \cdots & 0\\
-1 & -5 & 8 & -2 & \ddots & & & \vdots\\
0 & 0 & -6 & 8 & -2 & \ddots & & \vdots \\
\vdots & & \ddots & \ddots & \ddots & \ddots & \ddots & \vdots\\
\vdots & & & \ddots & \ddots & \ddots & \ddots & 0\\
\vdots & & & & \ddots & -6 & 8 & -2\\
\vdots & & & & & 0 & -5 & 5\\
0 & \cdots & \cdots & \cdots & \cdots & 0 & -1 & -1
\end{bmatrix},
\qquad
\mathbf{Q}_v^- = -(\mathbf{Q}_v^+)^{\!\top},
\end{equation*}
and
\begin{equation*}
\mathbf{D}_v^+ = \mathbf{H}_v^{-1}\!\left(\mathbf{Q}_v^+ + \tfrac{1}{2}\mathbf{B}_v\right),
\qquad
\mathbf{D}_v^- = \mathbf{H}_v^{-1}\!\left(\mathbf{Q}_v^- + \tfrac{1}{2}\mathbf{B}_v\right),
\end{equation*}
where
\begin{equation*}
\mathbf{H}_v = \Delta v\,\mathrm{diag}\!\bigl(\tfrac{1}{4},\tfrac{5}{4},1,\ldots,1,\tfrac{5}{4},\tfrac{1}{4}\bigr),
\qquad
\mathbf{B}_v = \mathrm{diag}(-1,0,\ldots,0,1).
\end{equation*}
These operators are second-order accurate in the interior and first-order accurate at the boundary, corresponding to an SBP(2,1) scheme, and satisfy the discrete SBP identity
\begin{equation*}
\mathbf{H}_v \mathbf{D}_v^+ + (\mathbf{H}_v \mathbf{D}_v^-)^{\!\top} = \mathbf{B}_v.
\end{equation*}
}

\end{appendices}

\bibliography{sn-bibliography}% common bib file

\end{document}